\documentclass[11pt,reqno]{amsart}
 
\newcommand{\mysection}[1]{\section{#1}
      \setcounter{equation}{0}}

\newcommand{\nlimsup}{\operatornamewithlimits{\overline{lim}}}

\newtheorem{theorem}{Theorem}[section]
\newtheorem{lemma}[theorem]{Lemma}

\newtheorem{corollary}[theorem]{Corollary}

\theoremstyle{definition}
\newtheorem{assumption}{Assumption}[section]
\newtheorem{definition}{Definition}[section]
\theoremstyle{remark}
\newtheorem{remark}{Remark}[section]

\newcommand\cbrk{\text{$]$\kern-.15em$]$}} 
\newcommand\opar{\text{\raise.2ex\hbox{${\scriptstyle | }$}\kern-.34em$($} }
\newcommand{\tr}{\text{\rm tr}\,}

 \makeatletter
 \def\dashint{%
 \operatorname%
 {\,\,\text{\bf--}\kern-.98em\DOTSI\intop\ilimits@\!\!}}
 \makeatother

\newcommand\gb{\mathfrak{b}}

\newcommand\bR{\mathbb{R}}

\newcommand\bL{\mathbb{L}}

\newcommand\bW{\mathbb{W}}

\newcommand\cB{\mathcal{B}}
\newcommand\cF{\mathcal{F}}

\newcommand\cP{\mathcal{P}}

\newcommand\cD{\mathcal{D}}

\newcommand\cL{\mathcal{L}}

\newcommand\cW{\mathcal{W}}

\newcommand\frD{\mathfrak{D}}

\begin{document}

\title[SPDEs
with growing   coefficients]{
On divergence form SPDEs
with growing  coefficients in $W^{1}_{2}$
spaces without weights}

\author[N.V.  Krylov]{N.V. Krylov}%
\thanks{The work   was partially supported
  by NSF grant DMS-0653121}
\address{127 Vincent Hall, University of Minnesota,
Minneapolis,
       MN, 55455, USA}
\email{krylov@math.umn.edu}


\begin{abstract}
We consider divergence form uniformly parabolic SPDEs with
bounded  and measurable  leading coefficients and
possibly growing lower-order coefficients in
the deterministic part of the equations. We look for solutions
  which are summable to the second power with respect to the usual
Lebesgue measure along with their first   derivatives
with respect to the spatial variable.
\end{abstract}

\maketitle

\mysection{Introduction}

We consider divergence form uniformly parabolic SPDEs with
bounded  and measurable  leading coefficients and
possibly growing lower-order coefficients
in the deterministic part of the equation.
 We look for solutions
  which are summable to the second power with respect to the usual
Lebesgue measure along with their first   derivatives
with respect to the spatial variable. To the best of our knowledge
our results are new even for deterministic PDEs when
one deletes all stochastic terms in the results below.
If there are no stochastic terms and the coefficients are 
nonrandom and
time independent, our results allow  one to
obtain the corresponding results for elliptic
divergence-form equations which also seem to be new.
A sample result in this case is the following.
Consider the equation
$$
D_{i}\big( a^{ij}( x)D_{j}u (x)
+\gb^{i} (x)u (x)\big) + b^{i} (x)D_{i} u (x)
$$
\begin{equation}
                                                     \label{6.9.1}
 -(c (x)+\lambda) u (x)=
D_{i}f^{i}(x)+f^{0}(x)
\end{equation}
in $\bR^{d}$ which is the Euclidean space
of points $x=(x^{1},...,x^{d})$. Here and below
 the summation convention
is enforced and
$$
D_{i}=\frac{\partial}{\partial x^{i}}.
$$
Assume that \eqref{6.9.1} is uniformly elliptic,
$a^{ij}$  are  bounded, and $c\geq0$. Also assume that
$f^{j}\in \cL_{2}=\cL_{2}(\bR^{d})$, $j=0,...,d$, and
$$
\sup_{|x-y|\leq1}(|b(x)-b(y)|+|\gb(x)-\gb(y)|+|c(x)-c(y)|)<\infty
$$
and that the constant $\lambda>0$ is large enough.
Then equation \eqref{6.9.1} has a unique solution in the class
of functions $u\in W^{1}_{2}=W^{1}_{2}(\bR^{d})$.
Notice that the above condition on $\gb,b$, and $c$
allow them to grow linearly as $|x|\to\infty$.

As in \cite{CV} one of the main motivations for studying
SPDEs with growing first-order coefficients
 is filtering theory 
for partially observable diffusion processes.
 
 It is generally believed that introducing weights is the 
most natural setting for equations with growing 
coefficients. When the 
coefficients grow it is quite natural to consider the equations
in function 
 spaces with weights that would restrict the set of solutions
in such a way that all terms in the equation will be
from the same space as the free terms. 
 The present paper seems to be the first
one treating the unique solvability of
these equations with growing
lower-order coefficients in the usual Sobolev spaces $W^{1}_{2}$
 without weights and
without imposing any {\em special\/}
 conditions on the relations between
the coefficients or on their {\em derivatives\/}.  

The theory of SPDEs in Sobolev-Hilbert spaces {\em
with\/} weights
attracted some attention in the past. We do not use weights and only
mention  a few papers about {\em stochastic\/}
PDEs in $\cL_{p}$-spaces
 with weights
in which one can find further references: \cite{AM} (mild solutions,
general $p$),
\cite{CV}, \cite{Gy93}, \cite{Gy97}, \cite{GK} ($p=2$ in the four
last articles).

Many more papers are devoted to the theory of {\em deterministic\/}
 PDEs
with growing coefficients in Sobolev spaces with weights. 
We  cite only a few of them sending the reader to the references
therein again because neither  do we deal with weights nor  use
the results of these papers. It is also worth saying that our results
do not generalize the results of the above cited papers.

In most of these papers the coefficients
are time independent, see \cite{CV1}, \cite{ChG}, \cite{FL}, \cite{Lu},
 \cite{MP1}, part of the result of which
are extended in \cite{GL} to
time-dependent Ornstein-Uhlenbeck operators.

It is worth noting that many issues for {\em deterministic\/}
divergence-type equations with time independent 
growing coefficients in $\cL_{p}$ spaces with arbitrary $p
\in(1,\infty)$ {\em without\/} weights
were also treated previously in the literature. 
This was done mostly  by using the semigroup approach
which excludes time dependent coefficients
and makes it almost impossible to use the results in
the more or less general filtering theory.
We briefly mention only a few recent papers
sending the reader to them for additional information.

In \cite{LV} a strongly continuous
in $\cL_{p}$ semigroup is constructed corresponding
to elliptic  operators with measurable
leading coefficients and Lipschitz
continuous drift coefficients.  
In \cite{MP} it is assumed that
 if, for $| x|\to\infty$, the drift coefficients   grow,
then  the zeroth-order coefficient  should grow, basically,
as the square of the drift. There  is also a condition on the divergence
of the drift coefficient.
In \cite{PR} there is no zeroth-order term
and the semigroup is constructed under some assumptions
one of which translates into  the monotonicity of 
$\pm b(x)-Kx$, for a constant $K$, if the leading term is the Laplacian.
In \cite{CF}  the drift coefficient
is assumed to be globally Lipschitz
continuous if the zeroth-order coefficient is constant.

Some conclusions in the above cited papers are quite similar to ours
but the corresponding assumptions are not as general
in what concerns the regularity of the coefficients.
However, these papers contain a lot of additional
important information not touched upon in the present paper
(in particular, it is shown in \cite{LV} that the corresponding semigroup
is not analytic).

The technique, we apply, originated from \cite{KP}
and \cite{Kr09_1} and uses special
 cut-off functions whose support evolves in time
in a manner adapted to the drift.
We do not make any regularity assumptions on the coefficients
and are  restricted to only treat equations in $W^{1}_{2}$. Similar,
techniques could be used to consider equations
in the spaces $W^{1}_{p}$ with any $p\geq2$. This time
one can use the results of \cite{Ki1} and \cite{Kr09_2} where  
some regularity on the coefficients in $x$ variable
is needed like, say, the condition that the 
second order coefficients
be in VMO uniformly with respect to the time variable
(see \cite{Kr09_2}).
However, for the sake
of brevity and clarity we concentrate only on 
$p=2$.  The main emphasis here is that we allow 
the first-order coefficients to grow
as $|x|\to\infty$ and still measure the size of the  
derivatives with respect to Lebesgue measure
thus avoiding using weights.
 
It is worth noting that considering
divergence form equations in $\cL_{p}$-spaces
is quite useful in the treatment of
filtering problems (see, for instance, \cite{Kr_10}) especially
when the power of summability is taken large  and we
intend to treat this issue in a subsequent paper.
 
The article is organized as follows. In Section
\ref{section 6.9.1} we describe the problem,
Section \ref{section 2.15.1} contains the statements
of two main results, Theorem \ref{theorem 3.11.1}
on an apriori estimate  providing,
in particular, uniqueness of solutions
 and Theorem \ref{theorem 3.16.1}
about the existence of solutions.
Theorem \ref{theorem 3.11.1} is proved in 
Section~\ref{section 6.9.3}
after we prepare the necessary tools in Section
\ref{section 6.9.2}.
Theorem \ref{theorem 3.16.1} is proved in the last
Section \ref{section 6.9.5}.

As usual when we speak of
 ``a constant" we always mean ``a finite constant".

\mysection{Setting of the problem}
                                       \label{section 6.9.1}

Let $(\Omega,\cF,P)$ be a complete probability space
with an increasing filtration $\{\cF_{t},t\geq0\}$
of complete with respect to $(\cF,P)$ $\sigma$-fields
$\cF_{t}\subset\cF$. Denote by $\cP$ the predictable
$\sigma$-field in $\Omega\times(0,\infty)$
associated with $\{\cF_{t}\}$. Let
 $w^{k}_{t}$, $k=1,2,...$, be independent one-dimensional
Wiener processes with respect to $\{\cF_{t}\}$.  
Finally, let $\tau$ be a stopping time with respect to
$\{\cF_{t}\}$.

 We consider the   second-order
operator $L$
\begin{equation}                                    \label{lu}
 L_{t} u_{t}(x) =D_{i}\big( a^{ij}_{t}( x)D_{j}u_{t}(x)
+\gb^{i}_{t}(x)u_{t}(x)\big) + b^{i} _{t}(x)
D_{i} u_{t}(x) -c _{t}(x) u_{t}(x),
\end{equation}
and the first-order operators
$$
\Lambda^{k}_{t} u_{t}(x)=\sigma^{ik}_{t}(x)D_{i}u_{t}(x)
+\nu^{k}_{t}(x)u_{t}(x)
$$  
   acting on functions $u_{t}(x)$ defined
on $\Omega\times\bR^{d+1}_{+}$, where $\bR^{d+1}_{+}=
 [0, \infty) \times \bR^d$, and given for $k=1,2,...$ 
(the summation convention is enforced throughout the article).
We set  $\bR_{+}=[0,\infty)$.

Our main concern is proving the unique solvability
of the equation
\begin{equation}
                                                \label{2.6.4}
du_{t}=(L_{t}u_{t}-\lambda u_{t}+D_{i}f^{i}_{t}+f^{0}_{t})\,dt
+(\Lambda^{k}_{t}u_{t}+g^{k}_{t})\,dw^{k}_{t},
\quad t\leq\tau,
\end{equation}
with an appropriate initial condition at $t=0$, where  
$\lambda>0$ is a constant. 
The precise 
assumptions on the coefficients, free terms, and initial data
will be given later. First we introduce appropriate function
spaces.

Denote $C^{\infty}_{0}=C^{\infty}_{0}(\bR^{d})$,
$\cL_{2}=\cL_{2}(\bR^{d})$, and let $W^{1}_{2} =W^{1}_{2}(\bR^{d})$
be the Sobolev space of functions $u$ of class
 $\cL_{2}$, such that
$Du\in \cL_{2}$, where $Du$ is the gradient of $u$. Introduce
$$
\bL_{2}( \tau)=\cL_{2}(\opar 0,\tau\cbrk,\bar{\cP},\cL_{2} ),
\quad
\bW^{1}_{2}( \tau)=\cL_{2}(\opar 0,\tau\cbrk,\bar{\cP},
W^{1}_{2} ),
$$
where $\bar{\cP}$ is the completion of $\cP$ with respect to the
product measure. 
Remember that the elements of
$\bL_{2}( \tau)$ need only  
belong to $\cL_{2}$ on a predictable subset of 
$\opar 0,\tau\cbrk$ of full measure. For the sake of convenience
we will always assume that they are defined everywhere
on $\opar 0,\tau\cbrk$ at least as generalized functions. 
Similar situation occurs in the case of $\bW^{1}_{2}( \tau)$.
We also use the same notation $\bL_{2}(\tau)$ for
$\ell_{2}$-valued functions like $g_{t}=(g^{k}_{t}
)$. For such a function, naturally,
$$
\|g\|_{\cL_{2}}=\|\,|g|_{\ell_{2}}
\,\|_{\cL_{2}}=\big\|\big(\sum_{k=1}^{\infty}
(g^{k})^{2}\big)^{1/2}\|_{\cL_{2}}
=\big(\sum_{k=1}^{\infty}\int_{\bR^{d}}
|g^{k}|^{2}\,dx\big)^{1/2}.
$$

  The following definition
turns out to be useful if the coefficients of
$L$ and $\Lambda^{k}$ are bounded. 
 
 \begin{definition}
                                         \label{definition 3.16.1}
     
We introduce the space $\cW^{1}_{2}(\tau)$,
which is the space of functions $u_{t}
=u_{t}(\omega,\cdot)$ on $\{(\omega,t):
0\leq t\leq\tau,t<\infty\}$ with values
in the space of generalized functions on $\bR^{d}$
and having the following properties:

(i) We have $u_{0}\in \cL_{2}(\Omega,\cF_{0},\cL_{2})$;

(ii)  We have $u
\in \bW^{1}_{2}(\tau )$;

(iii) There exist   $f^{i}\in \bL_{2}(\tau)$,
$i=0,...,d$, and $g=(g^{1},g^{2},...)\in \bL_{2}(\tau)$
such that
 for any $\phi\in C^{\infty}_{0}$ with probability 1
for all  $t\in\bR_{+}$
we have
$$
(u_{t\wedge\tau},\phi)=(u_{0},\phi)
+\sum_{k=1}^{\infty}\int_{0}^{t}I_{s\leq\tau}
(g^{k}_{s},\phi)\,dw^{k}_{s}
$$
\begin{equation}
                                                 \label{1.2.1}
+\int_{0}^{t}I_{s\leq\tau}\big(
 (f^{0}_{s},\phi)-(f^{i}_{s},D_{i}\phi)\big)\,ds.
\end{equation}
In particular, for any $\phi\in C^{\infty}_{0}$, the process
$(u_{t\wedge\tau},\phi)$ is $\cF_{t}$-adapted and (a.s.) continuous.
In case that property (iii) holds, we write
$$
du_{t}=(D_{i}f^{i}_{t}+f^{0}_{t})\,dt+g^{k}_{t}\,dw^{k}_{t},
\quad t\leq\tau.
$$
\end{definition}

It is a standard fact that for $g\in\bL_{2}(\tau)$
and any $\phi\in C^{\infty}_{0}$ the series
in \eqref{1.2.1} converges uniformly on $\bR_{+}$
in probability.

 Similarly  to this definition we understand
equation \eqref{2.6.4}  
in the general case  as the requirement that 
for any  $\phi\in C^{\infty}_{0}$ with probability one the relation
$$
(u_{t\wedge\tau},\phi)=(u_{0},\phi)
+\sum_{k=1}^{\infty}\int_{0}^{t}I_{s\leq\tau}
(\sigma^{ik}_{s}D_{i}u_{s}+\nu^{k}_{s}u_{s}
+g^{k}_{s},\phi)\,dw^{k}_{s}
$$
\begin{equation}
                                                \label{3.16.7}
+
\int_{0}^{t}I_{s\leq\tau}\big[(b^{i}_{s}D_{i}u_{s}
-(c_{s}+\lambda)u_{s}+f^{0}_{s},\phi)
-(a^{ij}_{s}D_{j}u_{s}+\gb^{i}_{s}u_{s}+
f^{i}_{s},D_{i}\phi)
 \big]\,ds
\end{equation}
hold   for all $t\in\bR_{+}$.

Observe that at this moment it is not clear that the right-hand
side makes sense.
Also notice that, if the coefficients of $L$ and $\Lambda^{k}$
are bounded, then any $u\in\cW^{1}_{2}(\tau)$ is a solution
of \eqref{2.6.4} with appropriate free terms since if
\eqref{1.2.1} holds, then \eqref{2.6.4} holds as well with
$$
f^{i}_{t}-a^{ij}_{t}D_{j}u_{t}-\gb^{i}u_{t},
\quad i=1,...,d,\quad
f^{0}_{t}+(c_{t}+\lambda)u_{t}-b^{i}_{t}D_{i}u_{t},
$$
$$
g^{k}_{t}-\sigma^{ik}D_{i}u_{t}-\nu^{k}_{t}u_{t}
$$
in place of $f^{i}_{t}$, $i=1,...,d$, $f^{0}_{t}$, and $g^{k}_{t}$,
respectively.

\mysection{Main results}
                                          \label{section 2.15.1}

For $\rho>0$ denote $B_{\rho}(x)=\{y\in\bR^{d}:|x-y|<\rho\}$,
$B_{\rho}=B_{\rho}(0)$. 
\begin{assumption} 
                                       \label{assumption 2.7.2}
(i) The functions $a^{ij}_{t}(x)$, $\gb^{i}_{t}(x)$,
$b^{i}_{t}(x)$, $c_{t}(x)$,
$\sigma^{ik}_{t}(x)$, $\nu^{k}_{t}(x)$ are real valued, measurable
with respect to $\cF\otimes\cB(\bR^{d+1}_{+})$,
 $\cF_{t}$-adapted for any $x$, and  $c\geq 0$.

(ii) There exist   constants $K,\delta>0$ such that
for all values of arguments and $\xi\in\bR^{d}$
$$
(a^{ij}   
-  \alpha^{ij} ) \xi^{i}
\xi^{j}\geq\delta|\xi|^{2},\quad
|a^{ij}|\leq \delta^{-1} ,
\quad  |\nu|_{\ell_{2}}\leq K,
$$
where  $\alpha^{ij} =(1/2)(\sigma^{i\cdot},\sigma^{j\cdot})
_{\ell_{2}}$. Also, the constant $\lambda>0$.

(iii) For any  $x\in
\bR^{d}$ (and $\omega$) the function
$$
\int_{B_{1}}(|\gb_{t}(x+y)|+|b_{t}(x+y)|+c_{t}(x+y))\,dy
$$ 
is locally square integrable on $\bR_{+}=[0,\infty)$.
\end{assumption}

Notice that the matrix $a=(a^{ij})$ need not be symmetric.
Also notice that in Assumption \ref{assumption 2.7.2} (iii)
the ball $B_{1}$ can be replaced with any other ball
without changing the set of admissible coefficients
$\gb,b,c$.

We take some 
$f^{j},g\in\bL_{2}(\tau)$ and before we give the definition of
solution of  \eqref{2.6.4}
we remind the reader
that, if $u\in\bW^{1}_{2}(\tau)$, then owing to the boundedness of $\nu$
and $\sigma$ and the fact that $Du,u,g\in\bL_{2}(\tau)$,
the first series on the right in \eqref{3.16.7}
converges uniformly in probability and the series
is a continuous local martingale.

\begin{definition}
                                      \label{definition 3.20.01}
By a solution of
\eqref{2.6.4} for $t\leq\tau$ with initial condition
$u_{0}\in\cL_{2}(\Omega,\cF_{0},\cL_{2})$
we mean a function $u\in \bW^{1}_{2}(\tau) $ 
 (not $\cW^{1}_{2}(\tau))$  such that

(i)  For any $\phi\in C^{\infty}_{0} $ with probability
one the integral
with respect to $ds$ in \eqref{3.16.7} is 
well defined and is finite for all
 $t\in\bR_{+}$;

(ii) For any $\phi\in C^{\infty}_{0} $
with probability one
equation  \eqref{3.16.7}  
 holds for all $t\in\bR_{+}$.
\end{definition}

For $d\ne2$ define
$$
q  =d\vee 2,
$$
and if $d=2$ let $q $ be a fixed number such that 
$q >2$.
The following assumption contains  a  parameter  $ 
\gamma \in(0,1]$,  
whose value will be specified later.

 \begin{assumption}[$\gamma $]
                                      \label{assumption 3.11.1} 
There exists a  $\rho_{0}\in(0,1]$ such that,
for any $\omega\in\Omega$
  and $ \gb:=(\gb^{1} ,...,\gb^{d} )$
and $ b:=(b^{1} ,...,b^{d} )$ and   $(t,x)\in\bR^{d+1}_{+}$
 we have  
$$
\rho_{0}^{- d}\int_{B_{\rho_{0}}}\int_{B_{\rho_{0}}}|\gb_{t}(  x+y)
-\gb_{t}(  x+z)|^{q }\,dydz  \leq\gamma ,
$$ 
$$
\rho_{0}^{- d}\int_{B_{\rho_{0}}}\int_{B_{\rho_{0}}}|b_{t}(  x+y)
-b_{t}(  x+z)|^{q}\,dydz  \leq\gamma ,
$$
$$
\rho_{0}^{- d}\int_{B_{\rho_{0}}}\int_{B_{\rho_{0}}}|c_{t}(  x+y)
-c_{t}(  x+z)|^{q}\,dydz  \leq\gamma .
$$
\end{assumption}

Obviously,  Assumption \ref{assumption 3.11.1}
is satisfied with any $\gamma\in(0,1]$ if  
 $b$, $\gb$, and $c$ are independent of $x$.
 Also notice that Assumption \ref{assumption 3.11.1}
 allows $b$, $\gb$, and $c$ growing linearly in $x$.

\begin{theorem}
                                       \label{theorem 3.11.1}
There exist
$$
\gamma =\gamma (d,\delta,K  )\in(0,1],
$$
$$
 N=N(d,\delta,K   ),
\quad \lambda_{0}=\lambda_{0}(d,\delta,K, 
\rho_{0})\geq0
$$  
such that, if the above assumptions are satisfied
and $\lambda\geq
\lambda_{0}$ and
  $u $ is a solution of \eqref{2.6.4}
with initial condition $u_{0}$  and some
$f^{j},g\in\bL_{2}(\tau)$, then
$$
\|u\sqrt{\lambda +c}\|^{2}_{\bL _{2}(\tau)}+\|Du\|^{2}_{\bL _{2}(\tau)}
\leq N\big(\sum_{i=1}^{d}\|f^{i}\|^{2}_{\bL _{2}(\tau)}
$$
\begin{equation}
                                       \label{3.11.2}
+\|g\|^{2}_{\bL _{2}(\tau)}
+\lambda^{-1}\|f^{0}\|^{2}_{\bL _{2}(\tau)}+
E\|u_{0}\|^{2}_{\cL_{2}}\big).
\end{equation}
\end{theorem}

This theorem provides an apriori estimate implying
uniqueness of solutions $u $.
Observe that the assumption that such a solution
exists is quite nontrivial because   if $\gb_{t}(x)\equiv x$,
it is not true that  $\gb u\in\bL_{2}(\tau)$
for arbitrary $u\in \bW^{ 1}_{2}(\tau) $.

To prove the existence we need stronger assumptions
because, generally, under only the above assumptions
the term
$$
D_{i}(\gb^{i}_{t}u_{t})+b^{i}_{t}D_{i}u_{t}
$$
cannot be written even locally as 
$D_{i}\hat{f}^{i}_{t}+\hat{f}^{0}_{t}$
with $\hat{f}^{j}\in\bL_{2}(\tau)$ if we only know that
$u\in\bW^{1}_{2}(\tau)$ even if $\gb$ and $b$ are independent
of $x$. 
We can only prove our crucial Lemma \ref{lemma 3.16.5}
 if such a representation is possible.

\begin{assumption}
                                       \label{assumption 3.16.1}
For any   $T,R\in\bR_{+}$, and $\omega\in\Omega$ we have  
$$
\sup_{t\leq T}\int_{B_{R}}(|\gb_{t}(x)| 
+|b_{t}(x)|+ c_{t}(x)   ) \,dx<\infty.
$$
 \end{assumption}

\begin{theorem}
                                        \label{theorem 3.16.1}
Let the above assumptions be satisfied with
$\gamma $ taken from Theorem \ref{theorem 3.11.1}.
Take   $\lambda\geq\lambda_{0}$, where $\lambda_{0}$
is defined in Theorem \ref{theorem 3.11.1}, and take 
$u_{0}\in\cL_{2}(\Omega,\cF_{0},\cL_{2})$. Then
there exists
a unique solution of \eqref{2.6.4}
  with initial condition $u_{0}$.

\end{theorem}

\begin{remark}
If the stopping time $\tau$ is bounded, then in the above theorem
one can take $\lambda=0$. To show this take a large $\lambda>0$
and replace
 the unknown function $u_{t}$ with $v_{t}
e^{ \lambda t}$. This leads to an equation for $v_{t}$
with the additional term $-\lambda v_{t}\,dt$ and the free terms
multiplied by $e^{-\lambda t}$. The existence
of $v\in\cW^{1}_{2}(\tau)$ will be then equivalent
to the existence of $u\in\cW^{1}_{2}(\tau)$ if $\tau$
is bounded.

\end{remark}

\mysection{A version of the It\^o-Wentzell formula}
                                             \label{section 6.9.2}

  Let $\cD$ be the space 
of
generalized functions on   $\bR^{d}$. 
We remind a definition and a result from \cite{Kr09_4}.
 Recall that for any $v\in\cD$
and $\phi\in C^{\infty}_{0}$ the function $(v,\phi(\cdot-x))$
is infinitely differentiable with respect to $x$, so that 
the sup in \eqref{11.16.2} below is predictable.

\begin{definition} 
 Denote by $\frD$                   \label{def 10.25.1}
 the set of all $\cD$-valued 
functions $u$ (written 
as $u_{t}(x)$ in a common abuse of notation)
on $\Omega\times\bR_{+}$ such that, for any $\phi\in C_{0}^{\infty}
:=C_{0}^{\infty}(\bR^{d})$,
 the restriction of the function $(u_{t},\phi)$ 
on $\Omega\times(0,\infty)$ is $\cP$-measurable 
and $(u_{0},\phi)$ is $\cF_{0}$-measurable.
 For $p=1,2$
denote by $\mathfrak{D}^{p}$ the subset of $\frD$
consisting of $u$ such that,  
  for any  
  $\phi\in C_{0}^{\infty}$  and
 $T ,R \in\bR_{+}$, we have
\begin{equation}
                                            \label{11.16.2}
\int_{0}^{T}\sup_{ |x|\leq R}|(u_{t} ,
\phi(\cdot-x))|^{ p}\,dt<\infty
\quad\hbox{ (a.s.)}.
\end{equation}
In the same way, considering $\ell_{2}$-valued
distributions $g$  on $C_{0}^{\infty}$, that is
linear $\ell_{2}$-valued functionals
such that $(g,\phi)$ is continuous as an $\ell_{2}$-valued
function with respect to the standard convergence of
test functions, we define 
$\frD(\ell_{2})$
 and  $\mathfrak{D}^{ 2}
(\ell_{2})$ 
replacing
$|\cdot|$ in (\ref{11.16.2})  
with $p=2$  by $|\cdot|_{\ell_{2}}$.

\end{definition}

Observe that if $g\in\mathfrak{D}^{2}(l_{2})$ then
for any $\phi\in C_{0}^{\infty}$,   and $T\in\bR_{+}$
$$
\sum_{k=1}^{\infty}\int_{0}^{T}(g^{k}_{t},\phi)^{2}\,dt
=\int_{0}^{T}|(g _{t},\phi)|_{\ell_{2}}^{2}\,dt<\infty
\quad\hbox{(a.s.)},
$$
which, by well known theorems about convergence of
series of martingales, implies that the series in \eqref{12.23.40} below
converges uniformly on $[0,T]$ in probability for any $T\in\bR_{+}$.

\begin{definition} 
                                           \label{def 10.25.3}
Let $f,u\in\mathfrak{D}$, $g\in\mathfrak{D}(l_{2})$.
 We say that the equality
\begin{equation}
                                           \label{11.16.3}
du_{t}(x)=f_{t}( x)\,dt+
g_{t}^{k}( x)\,dw^{k}_{t},\quad t\leq\tau,
\end{equation}
holds {\em in the sense of distributions\/} if 
$ fI_{\opar0,\tau\cbrk}\in\mathfrak{D}^ {1}$, 
$gI_{\opar0,\tau\cbrk}\in\mathfrak{D}^{ 2}(l_{2})$, and
for
 any $\phi\in C_{0}^{\infty}$, 
 with probability one  we have for all $t\in\bR_{+}$
\begin{equation}
                                             \label{12.23.40}
(u_{t\wedge\tau} ,\phi)=(u_{0} ,\phi)+\int_{0}^{t}I_{
s\leq\tau}
(f_{s},\phi)\,ds+\sum_{k=1}^{\infty}
\int_{0}^{t}I_{
s\leq\tau}(g^{k}_{s},\phi)\,dw^{k}_{s}.
\end{equation}
\end{definition}

Let $x_{t}$ be an $\bR^{d}$-valued stochastic process given by
$$
x^{i}_{t}=\int_{0}^{t}\hat b^{i}_{s}\,ds
+\sum_{k=1}^{\infty}\int_{0}^{t}\hat \sigma^{ik}_{s}\,dw^{k}_{s},
$$
where $\hat b_{t}=(\hat b^{i}_{t}),\hat \sigma^{k}_{t}
=(\hat \sigma^{ik}_{t})$ are predictable $\bR^{d}$-valued  
processes
such that for all $\omega$ and $s,T\in\bR_{+}$ we have
$\tr \hat\alpha_{s}<\infty$ and
$$
\int_{0}^{T}(|\hat b_{t}|+\tr \hat\alpha_{t})
\,dt<\infty,
$$
where $\hat\alpha_{t}=(\hat\alpha^{ij}_{t})$ and  
$2\hat \alpha^{ij}_{t}=(\hat \sigma^{i\cdot},
\hat \sigma^{j\cdot})_{\ell_{2}}$.
Finally, before stating the main result of \cite{Kr09_4}
we remind  the reader that for a generalized function $v$,
and any $\phi\in C^{\infty}_{0} $ the function
$(v,\phi(\cdot-x))$ is infinitely differentiable 
and for any derivative operator $D $ of order $n$ with respect 
to $x$ we have 
\begin{equation}
                                                 \label{4.3.2}
D (v,\phi(\cdot-x)) =(-1)^{n}(v,(D \phi)(\cdot-x))=:(D v,\phi(\cdot-x))=:
((D v)(\cdot+x),\phi)
\end{equation}
implying, in particular, that $D u \in\frD $ if $u \in\frD $.

\begin{theorem}
                                   \label{theorem 11.16.5}
Let $f,u\in\mathfrak{D}$, and $g\in\mathfrak{D}(l_{2})$.
Introduce
$$
v_{t}(x)=u_{t}(x+x_{t})
$$
and assume that \eqref{11.16.3} holds
 (in the sense of distributions). Then
$$
dv_{t}(x)= [f_{t}( x+x_{t})+\hat L_{t}v_{t}(x)
+(D_{i}g_{t}( x+x_{t}) ,\hat 
\sigma^{i\cdot}_{t} )_{\ell_{2}}]\,dt
$$
\begin{equation}
                                                  \label{4.4.5}
+
[g^{k}_{t}( x+x_{t})+D_{i}v_{t}( x)
\hat \sigma^{ik}_{t} ]\,dw_{t}^{k},
\quad t\leq\tau
\end{equation}
(in the sense of distributions), where $\hat L_{t}v_{t}=
\hat \alpha^{ij}_{t}
D_{i}D_{j}v_{t}( x)+\hat b^{i}_{t}D_{i}v_{t}(x)$.
In particular, the terms on the right in \eqref{4.4.5}
belong to the right class of functions.
\end{theorem}

We remind the reader  that the summation convention
over the repeated indices $i,j=1,...,d$ 
(and $k=1,2,...$) is enforced   throughout the
article.
  In the main part of this paper we are going to use
Theorem \ref{theorem 11.16.5} only for $ \hat\sigma \equiv0$.

\begin{corollary}
                                              \label{corollary 7.3.1}
Under the assumptions of Theorem \ref{theorem 11.16.5}
for any $\eta\in C^{\infty}_{0}$ we have
$$
d[u_{t}(x)\eta(x-x_{t})]=[g^{k}_{t}(x)\eta(x-x_{t})-
u_{t}(x)\hat\sigma^{ik}_{t}(D_{i}\eta)(x-x_{t})]\,dw^{k}_{t}
$$
$$
+[f_{t}(x)\eta(x-x_{t})+u_{t}(x)
(\hat{L}_{t}^{*}\eta)(x-x_{t})-(g_{t}(x),\hat\sigma^{i\cdot}(D_{i}\eta)
(x-x_{t}))_{\ell_{2}}]\,dt,\quad t\leq\tau,
$$
where $\hat{L}_{t}^{*}$ is the formal adjoint to $\hat{L}_{t}$.

\end{corollary}

Indeed, what we claim is that for any $\phi\in C^{\infty}_{0}$
with probability one
$$
((u_{t\wedge\tau} \phi)(\cdot+x_{t\wedge\tau}),\eta)=( u_{0} \phi ,\eta)
$$
$$
+\int_{0}^{t}I_{s\leq\tau} \big
(\big[ g^{k}_{s}\phi  +
\hat{\sigma}^{ik}_{s} D_{i}
(u_{s}\phi)\big](\cdot+x_{s}),\eta\big) \,dw^{k}_{s}
$$
$$
+\int_{0}^{t}I_{s\leq\tau}\big(\big[f_{s}\phi 
+ \hat{L}_{s}(u_{t}\phi)  +
 (\hat{\sigma}^{i\cdot}_{s},D_{i}(g _{s}\phi))_{\ell_{2}}\big]
(\cdot+x_{s}),\eta\big)\,ds
$$
for all $t$. However, to obtain this result it suffices
to write down an obvious equation for $u_{t}\phi$, then use 
Theorem \ref{theorem 11.16.5} and, finally,
use Definition \ref{def 10.25.3}
to interpret the result.

\mysection{Proof of Theorem \protect\ref{theorem 3.11.1}}
                                           \label{section 6.9.3} 

Throughout this section we suppose that the 
assumptions of Theorem~\ref{theorem 3.11.1} are satisfied and
 start with analyzing the second integral in
\eqref{3.16.7}. Recall that $q$ was introduced before
Assumption \ref{assumption 3.11.1}.

\begin{lemma}
                                        \label{lemma 6.27.1}
Let $h\in\cL_{q}$, $v\in
\cL_{2}$, and $u\in W^{1}_{2}$. Then   there exist
$V^{j}\in\cL_{2}$, $j=0,1,...,d$, such that
$$
h v=D_{i}V^{i}+V^{0},\quad
\sum_{j=0}^{d}\|V^{j}\|_{\cL_{2}}\leq
N\|h\|_{\cL_{q}}\|v\|_{\cL_{2}},
$$
where $N$ is independent of $h$ and $v$. In particular,
\begin{equation}
                                                         \label{7.1.1}
|(hv,u)|\leq N\|h\|_{\cL_{q}}\|v\|_{\cL_{2}}\|u\|_{W^{1}_{2}}.
\end{equation}
Furthermore, if   a number $\rho>0$,
then for any ball $B$ of radius $\rho$ we have
\begin{equation}
                                                 \label{6.27.7}
\|I_{B}hu\|_{\cL_{2}} \leq N\|h\| _{\cL_{q}}\big
(\rho^{  1-d/q }\|I_{B}Du\| _{\cL_{2}}
+\rho^{- d/q}\|I_{B}u\| _{\cL_{2}}\big),
\end{equation}
where $N$ is independent of $h$, $u$, $\rho$, and $B$.
\end{lemma}
Proof. Observe that by H\"older's inequality
for $r=2q/(2+q)$ ($\in[1,2)$) we have
$$
\|h v\|_{\cL_{r}}\leq  \|h\|_{\cL_{q}}\|v\|_{\cL_{2}}.
$$
Next we use the classical theory and introduce
$V \in W^{2}_{r}$
(note that $r>1$ if $d\ne1$ and $r=1$ if $d=1$)
as a unique solution of
$$
\Delta V - V =hv.
$$
 We know that
for a constant $N=N(d,r)$ we have
$$
\|V \|_{W^{2}_{r}}\leq N
\|h v\|_{\cL_{r}},
\quad
\|V \|_{W^{1}_{2}}\leq N\|V \|_{W^{2}_{r}},
$$
where the last inequality follows by embedding theorems
($2-d/r\geq1-d/2$). 
Now to prove the first assertion of the lemma
it only remains to combine the above estimates
and notice that for $V^{i}=D_{i}V$, $i=1,...,d$,
$V^{0}=-V$ it holds that $h v=D_{i}V^{i}+V^{0}$.

To prove the second assertion, first let $q>2$.
 Observe that by H\"older's inequality 
$$
\|I_{B}hu\|_{\cL_{2}}
\leq\|h\| _{\cL_{q}}\|I_{B}u\| _{\cL_{s}},
$$
where $s=2q/(q-2)$. By
embedding theorems (we use the fact that $d/s\geq d/2-1$)
$$
\|I_{B}u\| _{\cL_{s}}
\leq N(\rho^{ 1-d/q  }\|I_{B}Du\| _{\cL_{2}}
+\rho^{- d/q}\|I_{B}u\| _{\cL_{2}}\big)
$$
and the result follows. In the remaining case $q=2$,
which happens only if $d=1$. In that case the above
estimates remain  true if we set $s=\infty$.
The lemma is proved.

Before we extract some consequences from the lemma we 
take  a  nonnegative
$ \xi\in C^{\infty}_{0}(B_{\rho_{0}})$
  with unit integral and 
define
$$
\bar{b}_{s}(x)=\int_{B_{\rho_{0}}}\xi(y) b_{s}(x-y) \,dy,\quad
\bar{\gb}_{s}(x) =\int_{B_{\rho_{0}}}\xi(y) \gb_{s}(x-y) \,dy,
$$
\begin{equation}
                                                      \label{6.28.3}
\bar{c}_{s}(x)=\int_{B_{\rho_{0}}}\xi(y) c_{s}(x-y) \,dy.
\end{equation}
We may assume that $|\xi|\leq N(d)\rho_{0}^{-d}$.

One obtains the first two assertions of the following corollary
from \eqref{7.1.1} and \eqref{6.27.7}
by performing estimates like
$$
\|I_{B_{\rho_{0}}(x_{t})}(b_{t}-\bar{b}_{t}(x_{t}))\|_{\cL_{q}}^{q}
= \int_{B_{\rho_{0}}(x_{t})}| 
b_{ t}-\bar{b}_{t}(x_{t})|^{q}\,dx  
$$
$$
=\int_{B_{\rho_{0}}(x_{t})}\big|\int_{B_{\rho_{0}}(x_{t})}
[b_{t}(x)-b_{t}(y)]\xi(x_{t}-y)\,dy\big|^{q}\,dx
$$
$$
\leq N\int_{B_{\rho_{0}}(x_{t})}\big|\rho_{0}^{-d}
\int_{B_{\rho_{0}}(x_{t})}
|b_{t}(x)-b_{t}(y)|\,dy\big|^{q}\,dx
$$
\begin{equation}
                                                \label{6.11.2}
\leq N\rho_{0}^{-d}\int_{B_{\rho_{0}}(x_{t})} 
\int_{B_{\rho_{0}}(x_{t})}
|b _{t}(x)-b_{t}(y)|^{q}\,dy \,dx\leq N\gamma  ,
\end{equation}

\begin{corollary}
                                  \label{corollary 6.27.2}
Let $u\in\bW^{1}_{2}(\tau)$,
 let $x_{s}$ be
an $\bR^{d}$-valued predictable process, and let
$\eta\in C^{\infty}_{0}(B_{\rho_{0}})$. Set
$\eta_{s}(x)=\eta(x-x_{s})$. Then   on $\opar0,\tau\cbrk$

(i) For any $v\in W^{1}_{2}$
$$
 (|b^{i}_{s}-\bar{b}^{i}_{s}(x_{s})|I_{B_{\rho_{0}}(x_{s})}
| D_{i}u_{s}|,|v|) 
\leq N(d)\gamma^{1/q}\|I_{B_{\rho_{0}}(x_{s})}Du_{s}\| _{\cL_{2}}
\|v\|_{W^{1}_{2}} ;
$$

(ii) We have
$$
\|I_{B_{\rho_{0}}(x_{s})}
|\gb_{s}-\bar{\gb}_{s}(x_{s})|\,u_{s}\|_{\cL_{2}}
+\|I_{B_{\rho_{0}}(x_{s})}|c_{s}-\bar{c}_{s}(x_{s})|\,u_{s}\|_{\cL_{2}}
$$
$$
\leq N(d)\gamma^{1/q} \big
(\rho_{0}^{  1-d/q }\|I_{B_{\rho_{0}}(x_{s})}Du_{s}\| _{\cL_{2}}
+\rho_{0}^{- d/q}\|I_{B_{\rho_{0}}(x_{s})}u_{s}\| _{\cL_{2}}\big);
$$
 
(iii) Almost everywhere on $\opar0,\tau\cbrk$ we have
\begin{equation}
                                                \label{7.1.3}
(b^{i}_{s}-\bar{b}^{i}_{s}(x_{s}))\eta_{s} D_{i}u_{s}
=D_{i}V^{i}_{s}+V^{0}_{s} ,
\end{equation}
\begin{equation}
                                                \label{6.11.3}
\sum_{j=0}^{d}\|V^{j}_{s}\|_{\cL_{2}}
\leq N(d)\gamma^{1/q}\|I_{B_{\rho_{0}}}(x_{s})Du_{s}\|_{\cL_{2}}
\sup_{B_{\rho_{0}}}|\eta|,
\end{equation}
where $V^{j}_{s} $, $j=0,...,d$, are some predictable 
$\cL_{2}$-valued functions on $\opar0,\tau\cbrk$.

 \end{corollary}

To prove (iii) observe that one can find a predictable 
set $A\subset\opar0,\tau\cbrk$
of full measure such that  
$I_{A}D_{i}u$, $i=1,...,d$,
 are well defined as $\cL_{2}$-valued predictable
functions. Then \eqref{7.1.3} 
with $I_{A}D_{i}u$ in place of $D_{i}u$
and \eqref{6.11.3} follow  from \eqref{6.11.2},
 the first assertion of Lemma
\ref{lemma 6.27.1}, and the fact that the way $V^{j}$ are constructed
uses bounded hence continuous operators and  translates
the measurability of the data to the measurability of the result.
Since we are interested in
\eqref{7.1.3} and \eqref{6.11.3} holding only almost 
everywhere on $\opar0,\tau\cbrk$,
 there is no  actual need for the replacement.

\begin{corollary}
                                  \label{corollary 6.27.1}
Let $u\in\bW^{1}_{2}(\tau)$. Then 
for almost any $(\omega,s)$ the mappings
\begin{equation}
                                             \label{6.27.2}
\phi\,\to\,I_{s\leq\tau}(b^{i}_{s}D_{i}u_{s},\phi),\quad
I_{s\leq\tau}(\gb^{i}_{s} u_{s},D_{i}\phi),\quad
I_{s\leq\tau}(c_{s}u_{s},\phi)
\end{equation}
are generalized functions on $\bR^{d}$. Furthermore,
for any $T\in\bR_{+}$ almost surely
\begin{equation}
                                             \label{6.27.3}
\int_{0}^{T}
I_{s\leq\tau}(|(b^{i}_{s}D_{i}u_{s},\phi)|+
 |(\gb^{i}_{s} u_{s},D_{i}\phi)|+
|(c_{s}u_{s},\phi)|)\,ds<\infty,
\end{equation}
so that requirement (i) in Definition \ref{definition 3.20.01}
can be dropped.

\end{corollary}
 Proof. By having in mind partitions of unity we convince
ourselves that it suffices to prove that the mappings
\eqref{6.27.2} are generalized functions on any ball
$B$ of radius $\rho_{0}$ and that \eqref{6.27.3}
holds if $\phi\in C^{\infty}_{0}(B)$. Let $x_{0}$
be the center of $B$ and set $x_{s}\equiv x_{0}$. Then to prove the first assertion
concerning the last two functions in \eqref{6.27.2}
it suffices to use the first assertion of Corollary
\ref{corollary 6.27.2}  along with the observation that, say,
$$
(\gb^{i}_{s} u_{s},D_{i}\phi)=((\gb^{i}_{s}-\bar{\gb}^{i}_{s}(x_{0}))
 u_{s},D_{i}\phi)
+\bar{\gb}^{i}_{s}(x_{0})(  u_{s},D_{i}\phi).
$$
Similar transformation and   Corollary
\ref{corollary 6.27.2} (i) prove that the first function in 
\eqref{6.27.2} is also a generalized function.
  Assumption
\ref{assumption 2.7.2} (iii) and the estimates from Corollary
\ref{corollary 6.27.2} also easily imply \eqref{6.27.3}
thus finishing the proof of the corollary.

Before we continue with the proof of Theorem \ref{theorem 3.11.1},
 we notice that, if $u \in 
\cW^{ 1}_{2 }(\tau)$, then as we know (see, for instance,
 Theorem 2.1 of \cite{Kr09_3}),
there exists an event $\Omega'$ of full probability
such that $u_{t\wedge\tau}I_{\Omega'}$
is a continuous $\cL_{2}$-valued $\cF_{t}$-adapted process on
$\bR_{+}$. 
Substituting, $u_{t\wedge\tau}I_{\Omega'}$ in place of $u$ in our assumptions
and assertions does not change them. 
Furthermore, replacing $\tau$
with $\tau\wedge n$ and then sending $n$ to infinity
allows us to assume that $\tau$ is bounded.
Therefore, without losing
generality   we assume that

(H)   If we are considering a $u \in 
\cW^{ 1}_{2 }(\tau)$, the process
  $u_{t\wedge\tau}$ is a continuous $\cL_{2}$-valued 
$\cF_{t}$-adapted process on
$\bR_{+}$. The stopping time $\tau$ is bounded. 

Now we are ready to prove Theorem \ref{theorem 3.11.1}
in a particular case.

\begin{lemma}
                                        \label{lemma 3.11.1}
Let $ \nu^{k}\equiv0$ and let $\gb^{i}$, $b^{i}$,  and
$c$ be independent of $x$. Assume that
$u  $ is a solution of \eqref{2.6.4}
 with   some
$f^{j},g\in\bL_{2}(\tau)$  and $\lambda>0$.
 Then \eqref{3.11.2} holds with 
 $N=N(d,\delta,K)$.

\end{lemma}

Proof. 
We want to use Theorem \ref{theorem 11.16.5} to get rid
of the first order terms. Observe that \eqref{2.6.4} reads as
$$
du_{t}=(\sigma^{ik}_{t}D_{i}u_{t}+g^{k}_{t})\,dw^{k}_{t}  
$$
\begin{equation}
                                                \label{6.28.1}
+\big(D_{i}(a^{ij}_{t}D_{j}u_{t}+[\gb^{i}_{t} + 
b^{i}_{t}]u_{t}+f^{i}_{t})+f^{0}_{t}-(c_{t}
+\lambda) u_{t}\big)\,dt,
\quad t\leq\tau.
\end{equation}

One can find a predictable set $A\subset\opar0,\tau\cbrk$
of full measure such that $I_{A}f^{j}$, $j=0,1,...,d$, and
$I_{A}D_{i}u$, $i=1,...,d$,
 are well defined as $\cL_{2}$-valued predictable
functions satisfying
$$
\int_{0}^{\infty}I_{A}\big(\sum_{j=0}^{d}\|f^{j}_{t}\|^{2}_{\cL_{2}}
+ \|Du_{t}\|^{2}_{\cL_{2}}\big)\,dt<\infty.
$$
Replacing $f^{j}$ and $D_{i}u$ in \eqref{6.28.1} with
$I_{A}f^{j}$ and $I_{A}D_{i}u$, respectively, will not affect
\eqref{6.28.1}. Similarly, one can handle the function $g$
  and the terms $h_{t}=I_{\opar0,\tau\cbrk}[\gb^{i}  + 
b^{i} ]u ,I_{\opar0,\tau\cbrk}c u $ for which
$$
\int_{0}^{T}\|h_{t}\|_{\cL_{1}}\,dt<\infty\quad\text{(a.s.)}
$$
for each $T\in\bR^{d}$ owing to Assumption
\ref{assumption 2.7.2} (iii) and the fact that
$u\in\bW^{1}_{2}(\tau)$.
 
After these replacements all terms in \eqref{6.28.1}
will be of   class $\frD^{1}$ or $\frD^{2}(\ell_{2})$
as appropriate since  
 $a$ and $\sigma$ are bounded. 
This allows us to apply Theorem \ref{theorem 11.16.5}
and for
$$
B_{t}^{i}=\int_{0}^{t}
(\gb^{i}_{s}+b^{i}_{s})\,ds,\quad 
 \hat{u}_{t}(x)=u_{t}(x-B_{t})
$$
obtain that
$$
d\hat{u}_{t}=\big(D_{i}(\hat{a}^{ij}_{t}D_{j}\hat{u}_{t} )
 -(c_{t}+\lambda) \hat{u}_{t}+D_{i}\hat{f}^{i}_{t}+\hat{f}^{0}_{t}\big)\,dt
$$
\begin{equation}
                                               \label{6.28.8}
+\big(\hat{\sigma}^{ik}_{t}D_{i}\hat{u}_{t}+
\hat{g}^{k}_{t}\big)\,dw^{k}_{t},\quad t\leq\tau,
\end{equation}
where
$$
(\hat{a}^{ij}_{t},\hat{\sigma}^{ik}_{t}, 
\hat{f}^{j}_{t},
\hat{g}^{k}_{t})(x)=(a^{ij}_{t},\sigma^{ik}_{t}, f^{j}_{t},
g^{k}_{t})(x-B_{t}).
$$

Obviously, $\hat{u}$ is in $\bW^{1}_{2}(\tau)$ and its norm
coincides with that of $u$.
  Moreover, having in mind that $c_{t}$
is independent of $x$ and is locally (square) integrable,
one can find stopping times
$\tau_{n}\uparrow\tau$ such that $I_{\tau_{n}\ne\tau}
\downarrow0$ and
$$
\xi_{\tau_{n}}\leq n,\quad \xi_{t}:=\int_{0}^{t}c_{s}\,ds\leq n .
$$
Then it follows from from the equation
$$
d(\xi_{t}\hat{u}_{t})=\big(D_{i}(
\xi_{t}\hat{a}^{ij}_{t}D_{j}\hat{u}_{t} )
 -\lambda \xi_{t}\hat{u}_{t}+D_{i}\xi_{t}\hat{f}^{i}_{t}+
\xi_{t}\hat{f}^{0}_{t}\big)\,dt
$$
$$
+\big(\hat{\sigma}^{ik}_{t}\xi_{t}D_{i}\hat{u}_{t}+
\xi_{t}\hat{g}^{k}_{t}\big)\,dw^{k}_{t},\quad t\leq\tau_{n}
$$
that $\xi u\in\cW^{1}_{2}(\tau_{n})$ and hence
$\xi_{t\wedge\tau_{n}} u _{t\wedge\tau_{n}}$
is a continuous $\cL_{2}$-valued function and so are
$u_{t\wedge\tau_{n}}$ and $u_{t\wedge\tau}$.

 Furthermore, 
 since $\tau$ is bounded and $u_{t\wedge\tau}$
is a continuous $\cL_{2}$-valued function 
 and $c_{t}$ is independent of
$x$ and is locally square integrable, we have
\begin{equation}
                                                     \label{6.29.2}
\int_{0}^{\tau}\|c_{t}\hat{u}_{t}\|^{2}_{\cL_{2}}\,dt
=\int_{0}^{\tau}c^{2}_{t}\|u_{t}\|^{2}_{\cL_{2}}\,dt
\leq\sup_{t\leq\tau}\|u_{t}\|^{2}_{\cL_{2}}\int_{0}^{\tau}
 c^{2}_{t}\,dt<\infty
\end{equation}
and there is a sequence of,  
perhaps, different from the above 
 stopping times $\tau_{n}\uparrow\tau$
such that for each $n$
\begin{equation}
                                                     \label{6.29.1}
E\int_{0}^{\tau_{n}}\|c_{t}\hat{u}_{t}\|^{2}_{\cL_{2}}\,dt<\infty.
\end{equation}
Then \eqref{6.28.8} implies that $\hat{u}\in\cW^{1}_{2}(\tau_{n})$
for each $n$. Also observe that if we can prove \eqref{3.11.2}
with $\tau_{n}$ in place of $\tau$, then we can
let $n\to\infty$ and use the monotone convergence theorem
to get \eqref{3.11.2} as is. Therefore, in the rest of the proof
we assume that \eqref{6.29.1} holds with $\tau$ in place of
$\tau_{n}$, that is,  $\hat{u}\in\cW^{1}_{2}(\tau )$.

The next argument is standard (see, for instance,  
 Lemma 3.3 and Corollary 3.2 of \cite{Kr09_2}). It\^o's formula
implies that
\begin{equation}
                                                  \label{3.11.3}
E\|u_{0}\|^{2}_{\cL_{2}}+
E\int_{0}^{\tau}\int_{\bR^{d}}I_{t}\,dxdt\geq0,
\end{equation}
where
$$
I_{t}:=2\hat{u}_{t}(\hat{f}^{0}_{t}-\lambda \hat{u}_{t}-c_{t}
\hat{u}_{t})
-2(\hat{a}^{ij}_{t}D_{j}\hat{u}_{t}+\hat{f}^{i}_{t})D_{i}\hat{u}_{t}
+|\hat{\sigma}^{i\cdot}_{t}D_{i}\hat{u}_{t}+
\hat{g}_{t}|_{\ell_{2}}^{2}.
$$
We use the inequality
$$
|\hat{\sigma}^{i\cdot}_{t}D_{i}\hat{u}_{t}+
\hat{g}_{t}|_{\ell_{2}}^{2}\leq (1+\varepsilon)
 | \hat\sigma^{i\cdot}_{t}D_{i}\hat{u}_{t} |^{2}
_{\ell_{2}}
+2\varepsilon^{-1}|\hat{g}_{t}|^{2}_{\ell_{2}},
\quad\varepsilon\in(0,1],
$$
and Assumption \ref{assumption 2.7.2}. Then for 
$\varepsilon=\varepsilon(\delta)>0$ small enough we find
$$
I_{t}\leq-\delta|D\hat{u}_{t}|^{2}-2(c_{t}+\lambda)\hat{u}^{2}_{t}
+2\hat{u}_{t}\hat{f}^{0}_{t}-2\hat{f}^{i}_{t} D_{i}\hat{u}_{t}
+N|\hat{g}_{t}|^{2}_{\ell_{2}}.
$$
Once again using $2\hat{u}_{t}\hat{f}^{0}_{t}\leq\lambda
\hat{u}^{2}_{t}+\lambda^{-1}|\hat{f}^{0}_{t}|^{2}$ and
similarly estimating $2\hat{f}^{i}_{t} D_{i}\hat{u}_{t}$
we conclude that
$$
I_{t}\leq-(\delta/2)|D\hat{u}_{t}|^{2}-(c_{t}+
 \lambda)\hat{u}^{2}_{t}
+N\big(\sum_{i=1}^{d}|\hat{f}^{i}_{t}|^{2}+|\hat{g}_{t}|
_{\ell_{2}}^{2}\big)+N\lambda^{-1}|\hat{f}^{0}_{t}|^{2}.
$$
By coming back to \eqref{3.11.3} we obtain
$$
 \|\hat{u}\sqrt{c_{t}+
 \lambda}\|^{2}_{\bL_{2}(\tau)}+\|D\hat{u}\|^{2}
_{\bL_{2}(\tau)}\leq N
 \big(\sum_{i=1}^{d}
\|\hat{f}^{i}\|^{2}_{\bL_{2}(\tau)}+\|\hat{g}\|^{2}_{\bL_{2}(\tau)}
\big)
$$
$$
+N\lambda^{-1 }\|\hat{f}^{0}\|^{2}_{\bL_{2}(\tau)}+N
E\|u_{0}\|^{2}_{\cL_{2}}.
$$
This is equivalent to \eqref{3.11.2} and the lemma is proved.

To proceed further we need a construction. Take $\bar{\gb},
\bar{b}$, and $\bar{c}$
from \eqref{6.28.3}.
 From Lemma 4.2 of \cite{Kr09_1} and Assumption
\ref{assumption 3.11.1} it follows  that, for 
$h_{ t}=\bar{\gb}_{ t},\bar{b}_{  t},
\bar{c}_{ t}$, it holds that  
 $|D^{n}h_{ t}|\leq \kappa_{n} $, where
$\kappa_{n}=\kappa_{n}
(n,\gamma ,d,\rho_{0})\geq1$ and $ D^{n}h_{ t }$
is any derivative of $h_{ t}$ of order $n\geq1$
with respect to $x$. By Corollary 4.3 of \cite{Kr09_1}
we have $|h_{ t}( x )|\leq K(t)(1+|x|)$,
where for each $\omega$ the function $K(t)=K(\omega,t)$
is locally 
(square) integrable with respect to $t$ on $\bR_{+}$.
Owing to these properties
the equation
\begin{equation}
                                             \label{2.8.1}
x_{t}=x_{0}-\int_{t_{0}}^{t}(\bar{\gb}_{ s}+\bar{b}_{  s})
( x_{s})\,ds,\quad t \geq t_{0} ,
\end{equation}
 for any ($\omega$ and)  
$ (t_{0},x_{0}) \in \bR^{d+1 }_{+} $  has a unique solution 
$x_{t}=x_{t_{0},x_{0},t} $.
Obviously, the process $x_{t_{0},x_{0},t}$, $t\geq t_{0}$,
 is $\cF_{t}$-adapted.

Next, for $i=1,2$ set $\chi^{(i)}(x)$ to be the indicator function
of $B_{\rho_{0}/i}$ and introduce
$$
\chi^{(i)}_{t_{0},x_{0},t}(x)=\chi^{(i)}(x-x_{t_{0},x_{0},t})
I_{t\geq t_{0}}.
$$

Here is a crucial estimate.
\begin{lemma}
                                     \label{lemma 3.14.1}
Assume that
$u $ is a solution of \eqref{2.6.4}
  with   some
$f^{j},g\in\bL_{2}(\tau)$. Then for $ (t_{0},x_{0}) \in \bR^{d+1 }_{+} $
and $\lambda>0$ we have
$$
  \|\chi^{(2)}_{t_{0},x_{0}}u\sqrt{c+\lambda}\|^{2}_{\bL _{2}(\tau)}+
\|\chi^{(2)}_{t_{0},x_{0}}Du\|^{2}_{\bL _{2}(\tau)}
$$
$$
\leq N\big(\sum_{i=1}^{d}\|\chi^{(1)}_{t_{0},x_{0}}
f^{i}\|^{2}_{\bL _{2}(\tau)}
+\|\chi^{(1)}_{t_{0},x_{0}}g\|^{2}_{\bL _{2}(\tau)}\big)
$$
$$
+N\lambda^{-1 }\|\chi^{(1)}_{t_{0},x_{0}}f^{0}\|^{2}_{\bL _{2}(\tau)}
+NE\|u_{t_{0}}I_{B_{\rho_{0}}(x_{0})}
I_{t_{0}\leq\tau}\|^{2}
_{\cL_{2}}
$$
$$
 +N\gamma ^{2/q} \| 
\chi^{(1)}_{t_{0},x_{0}} Du
\|_{\bL_{2}(\tau)}^{2}+
 N^{*} \lambda^{-1}\| 
\chi^{(1)}_{t_{0},x_{0}} Du
\|_{\bL_{2}(\tau)}^{2}
$$
\begin{equation}
                                       \label{3.14.2}
+ N^{*}(1+\lambda^{-1}) \|
\chi^{(1)}_{t_{0},x_{0}}  u
\|_{\bL_{2}(\tau)}^{2}
+N^{*}\lambda^{-1}\sum_{i=1}^{d}\|\chi^{(1)}_{t_{0},x_{0}}
f^{i}\|^{2}_{\bL _{2}(\tau)},
\end{equation}
where and below  in the proof
 by $N$ we denote generic constants depending only
on $d,\delta$, and $K$ and by $N^{*}$ constants depending only
on the same objects and $\rho_{0}$.
\end{lemma}

Proof. Since we are only concerned with
the values of $u_{t}$ if $t_{0}\leq t\leq\tau$,
 we may start considering
\eqref{2.6.4} on $[t_{0},\tau\vee t_{0})$ and then shifting time
allows us to assume that $t_{0}=0$. Obviously, we may also assume that
$x_{0}=0$. With this stipulations we will drop the subscripts $t_{0},
x_{0}$. Then, we can include the term    $\nu^{k}u$ into
  $g^{k}$  and obtain  \eqref{3.14.2}
by the triangle inequality
if we assume that this estimate is true in case 
$ \nu^{k}\equiv0$.
Thus, without losing generality we assume
$$
t_{0}=0,\quad x_{0}=0,\quad
 \nu^{k}\equiv0.
$$

Fix  a 
$ \zeta\in C^{\infty}_{0} $ with support in $B_{\rho_{0}}$
and such that $\zeta =1$ on $B_{\rho_{0}/2}$ and $0\leq\zeta
\leq1$. 
Set $x_{t}=x_{0,0,t}$,
$$
  \hat{\gb}_{ t} =\bar{\gb}_{ t}( x_{ t}) ,\quad
\hat{b}_{ t} =\bar{b}_{  t}( x_{ t}),\quad
\hat{c}_{ t} =\bar{c}_{  t}( x_{ t})
$$
$$
\eta_{ t}( x)=\zeta(x-x_{ t} ), 
\quad v_{ t}( x)= u_{t}( x)  \eta_{ t}( x).
$$
The most important property of $\eta_{t}$ is that
$$
d\eta_{t}=(\hat{\gb}^{i}_{t}+\hat{b}^{i}_{t})D_{i}\eta_{t}\,dt.
$$
Also
observe  for the later that we may assume that
\begin{equation}
                                                   \label{6.16.1}
\chi^{(2)}_{t}\leq\eta_{ t}\leq \chi^{(1)}_{t},\quad
|D\eta_{ t}|\leq N\rho_{0}^{-1 }\chi^{(1)}_{t},
\end{equation}
where $\chi^{(i)}_{t}=\chi^{(i)}_{0,0,t}$ and $N=N(d)$.

By Corollary \ref{corollary 7.3.1} (also see the argument before
\eqref{6.28.8})
we obtain
that  for $t\leq\tau$
$$
dv_{ t} =\big[D_{i}(\eta_{ t}a^{ij}_{t}D_{j}u_{t}+\gb^{i}_{t}v_{ t})
-(a^{ij}_{t}D_{j}u_{t}+\gb^{i}_{t}u_{t})
D_{i}\eta_{ t}
$$
$$
+b^{i}_{t}\eta_{ t}D_{i} u_{t} 
-(c_{t}+\lambda) v_{ t}
+D_{i}(f^{i}_{t}\eta_{ t})-f^{i}_{t}D_{i}\eta_{ t}
+f^{0}_{t}\eta_{ t}
$$
$$
+
(\hat{\gb}^{i}_{ t} +\hat{b}^{i}_{ t} )u_{t}
D_{i} \eta_{ t}\big]\,dt+\big[\sigma^{ik}D_{i}v_{ t}
-\sigma^{ik}u_{t}D_{i}\eta_{ t}+g^{k}_{t}\eta_{ t}
\big]\,dw^{k}_{t}.
$$
We transform  this further by noticing that  
$$
\eta_{ t}a^{ij}_{t}D_{j}u_{t}= a^{ij}_{t}D_{j}v_{ t}-
a^{ij}_{t}u_{t}D_{j}\eta_{ t}.
$$

To deal with the term $b^{i}_{t}\eta_{ t}D_{i} u_{t} $
we use Corollary \ref{corollary 6.27.2} and find 
the corresponding functions $V^{j}_{t}$.
Then simple arithmetics show that
$$
dv_{ t}=(\sigma^{ik}D_{i}v_{ t}
 +\hat{g}^{k}_{t} 
)\,dw^{k}_{t}
$$
$$
+\big[D_{i}\big(a^{ij}_{t}D_{j}v_{ t}+
\hat{\gb}^{i}_{ t}v_{ t}\big)-(\hat{c}_{ t}+\lambda) v_{ t}
+\hat{b}^{i}_{ t}
D_{i}v_{ t}+D_{i}\hat{f}^{i}_{ t}
+\hat{f}^{0}_{ t}\big]\,dt,
$$  
where
$$
\hat{f}^{0}_{ t}=f^{0}_{t}\eta_{ t}-f^{i}_{t}D_{i}\eta_{ t}
 -a^{ij}_{t}(D_{j}u_{t})D_{i}\eta_{ t}
+( \hat{\gb}^{i}_{ t}-\gb^{i}_{t})
u_{t}D_{i} \eta_{ t}+(\hat{c}_{ t}-c_{t})u_{t}\eta_{t}
+V^{0}_{t},
$$ 
$$
\hat{f}^{i}_{ t}=f^{i}_{t}\eta_{ t}-
a^{ij}_{t}u_{t}D_{j}\eta_{ t}+(\gb^{i}_{t}-
\hat{\gb}^{i}_{ t})
u_{ t}\eta_{ t}+V^{i}_{t},\quad i=1 ,..,d,
$$
$$
\hat{g}^{k}_{ t}=
-\sigma^{ik}u_{t}D_{i}\eta_{ t}+g^{k}_{t}\eta_{ t}.
$$

It follows by Lemma \ref{lemma 3.11.1} that for $\lambda>0$
$$
 \|v \sqrt{\hat{c}+\lambda}\|^{2}_{\bL_{2}(\tau)}+
\|Dv \|^{2}_{\bL_{2}(\tau)}\leq 
 N\lambda^{-1}\|\hat{f}^{0} \|^{2}_{\bL_{2}(\tau)}
$$
\begin{equation}
                                                   \label{3.13.1}
+N\big(
\sum_{i=1}^{d}\|\hat{f}^{i} \|^{2}_{\bL_{2}(\tau)}
+\|\hat{g} \|^{2}_{\bL_{2}(\tau)}+E\|v_0\|^{2}_{\cL_{2}}\big).
\end{equation}
Recall that here
  and below by $N$ we denote generic constants
depending only on $d,\delta$, and $K$.

Now we start estimating the right-hand side of \eqref{3.13.1}.
First we deal with $\hat{f}^{i}_{ t}$ and 
$\hat{g}^{k}_{ t}$.
Recall \eqref{6.16.1} and observe that obviously, if  
$\eta_{ t} (x)\ne0$, then $|x -x_{ t}|\leq\rho_{0}$.
Therefore,
\begin{equation}
                                          \label{3.14.02}
\|\hat{g}  \|^{2}_{\bL_{2}(\tau)}
\leq N^{*} \|u \chi^{(1)}_{\cdot}\|^{2}_{\bL_{2}(\tau)}
+N\|g \chi^{(1)}_{\cdot}\|^{2}_{\bL_{2}(\tau)}  
\end{equation}
(we remind the reader that
 by $N^{*}$ we denote generic constants depending
only on $d,\delta, K$, and $\rho_{0}$).
By Corollary \ref{corollary 6.27.2}  
\begin{equation}
                                          \label{3.14.1}
  \| (\gb^{i}_{t}-
\hat{\gb}^{i}_{ t})
u_{ t}\eta_{ t}\|^{2}_{\cL_{2}}
\leq N \gamma^{2/q}(\rho_{0}^{2(1-d/q)}\|\chi^{(1)}_{t}
 Du_{t}  \|^{2}_{\cL_{2}}+
 \rho_{0}^{-2d/q}\|\chi^{(1)}_{t}
  u_{t}\|^{2}_{\cL_{2}}).
\end{equation}  
Here $\rho_{0}^{2(1-d/q)}\leq1$ since $q\geq d$.
By adding that
$$
\|a^{ij} u D_{j}\eta \|^{2}
_{\bL_{2}(\tau)}\leq N^{*}
\|\chi^{(1)}_{\cdot}u \|^{2}_{\bL_{2}(\tau)},
$$
we derive from \eqref{6.11.3}, \eqref{3.14.02}, and \eqref{3.14.1}
that  
$$
\sum_{i=1}^{d}\|\hat{f}^{i} \|^{2}_{\bL_{2}(\tau)}
+\|\hat{g}  \|^{2}_{\bL_{2}(\tau)}
\leq N\big(\sum_{i=1}^{d}\|\chi^{(1)}_{\cdot}f^{i}  
\|^{2}_{\bL_{2}(\tau)}
+\|\chi^{(1)}_{\cdot}g  \|^{2}_{\bL_{2}(\tau)}\big)
$$
\begin{equation}
                                          \label{3.14.4}
+N\gamma ^{2/q} \| 
\chi^{(1)}_{\cdot}  Du
\|_{\bL_{2}(\tau)}^{2}+ N^{*}\|
\chi^{(1)}_{\cdot}  u
\|_{\bL_{2}(\tau)}^{2}.
\end{equation}

While estimating $\hat{f}^{0}$ we 
use \eqref{6.11.3} again and
observe that  we can deal with $( \hat{\gb}^{i}_{ t}-\gb^{i}_{t})
u_{t}D_{i} \eta_{ t}$ as in \eqref{3.14.1} this time without paying
much attention to the dependence of our constants on $\rho_{0}$
and obtain that
$$
\|( \hat{\gb}^{i} -\gb^{i} )
u D_{i} \eta \|_{\bL_{2}(\tau)}^{2}
\leq  N^{*}(\|\chi^{(1)}_{\cdot} Du
\|_{\bL_{2}(\tau)}^{2}+\|\chi^{(1)}_{\cdot}u\|_{\bL_{2}(\tau)}^{2}).
$$
By estimating also roughly the remaining terms in $\hat{f}^{0}$
and combining this with \eqref{3.14.4} and \eqref{3.13.1},
we see that the left-hand side of \eqref{3.13.1}
is less than the right-hand side of \eqref{3.14.2}.
  However,
$$
|\chi^{(2)}_{t}Du_{t}|\leq|\eta_{t}Du_{t}|\leq
|Dv_{t}|+|u_{t}D\eta_{t}|\leq|Dv_{t}|+
N\rho_{0}^{-1}|u_{t} \chi^{(1)}_{t}|
$$
and also
$$
|\chi^{(2)}_{t}u_{t}|^{2}(c_{t}+\lambda)\leq
|\eta_{t} u_{t}|^{2}(c_{t}+\lambda)
\leq |v_{t}|^{2}(\hat{c}_{t}+\lambda)+|\eta_{t} u_{t}|^{2}(1+|c_{t}
-\hat{c}_{t}|^{2}).
$$
By combining this with
the fact that by Corollary \ref{corollary 6.27.2}
$$
\|( \hat{c}^{i} -c )
u \eta \|_{\bL_{2}(\tau)}^{2}
\leq  N\gamma^{2/q}\|\chi^{(1)}_{\cdot} Du
\|_{\bL_{2}(\tau)}^{2}+N^{*} \|\chi^{(1)}_{\cdot}   u\|_{\bL_{2}(\tau)}^{2}) 
$$
we obtain \eqref{3.14.2}.
The lemma is proved.

Next, from  the  result giving ``local" in space estimates
we derive global in space estimates but for functions
having, roughly
speaking, small ``future" support in the time variable.  

\begin{lemma}
                                                \label{lemma 3.14.3}

Assume that
$u $ is a solution of \eqref{2.6.4}
  with   some
$f^{j},g\in\bL_{2}(\tau)$ and assume that $u_{t}=0$ if $t_{0}+\kappa_{1}^{-1}
\leq t\leq\tau$ with
$\kappa_{1}=\kappa_{1}(\gamma ,d,\rho_{0})\geq1$ introduced 
before \eqref{2.8.1} and some (nonrandom) $t_{0}\geq 0$
(nothing is required for those $\omega$ for 
which $\tau<t_{0}+\kappa^{-1}$).
 Then for $\lambda>0$ and $I_{t_{0}}:=I_{[t_{0},\infty)}$
$$
\| I_{t_{0}}u\sqrt{c+\lambda }\|^{2}_{\bL _{2}(\tau)}+
\|I_{t_{0}} Du\|^{2}_{\bL _{2}(\tau)}
\leq N\big(\sum_{i=1}^{d}\| 
I_{t_{0}}f^{i}\|^{2}_{\bL _{2}(\tau)}
+\|I_{t_{0}} g\|^{2}_{\bL _{2}(\tau)}\big)
$$
$$
+N\lambda^{-1 }\|I_{t_{0}} f^{0}\|^{2}_{\bL _{2}(\tau)}
+N
E\|u_{t_{0}}I_{t_{0}\leq\tau}\|^{2}_{\cL_{2}}
$$
$$
 +N\gamma ^{2/q} \| I_{t_{0}} Du\|_{\bL_{2}(\tau)}^{2}
+
 N^{*} \lambda^{-1}\| I_{t_{0}}  Du\|_{\bL_{2}(\tau)}^{2}
$$
\begin{equation}
                                       \label{3.14.5} 
+ N^{*}(1+\lambda^{-1})\| I_{t_{0}} u\|_{\bL_{2}(\tau)}^{2}
+N^{*}\lambda^{-1}\sum_{i=1}^{d}\| 
I_{t_{0}}f^{i}\|^{2}_{\bL _{2}(\tau)},
\end{equation}
where and below in the proof by $N$ we denote generic constants depending only
on $d,\delta$, and $K$ and by $N^{*}$ constants depending only
on the same objects and $\rho_{0}$.

\end{lemma}

Proof. Take $x_{0}\in\bR^{d}$ and use the notation introduced before
Lemma \ref{lemma 3.14.1}.
One knows that for each   $t\geq t_{0}$,
the mapping $x_{0}\to x_{t_{0},x_{0},t} $ is a diffeomorphism
with Jacobian determinant given by
$$
\bigg|\frac{\partial x_{t_{0},x_{0},t} }{
\partial x_{0}}\bigg| =\exp\big(-\int_{t_{0}}^{t}\sum_{i=1}^{d} D_{i}
[\bar{\gb} _{  s}^{i}+\bar{ b} _{  s}^{i}]
( x_{t_{0},x_{0},s}) \,ds\big).
$$
By the  way the constant $\kappa_{1}$ is introduced, we have
$$
e^{-N\kappa_{1}(t-t_{0})}\leq \bigg|\frac{\partial x_{t_{0},x_{0},t}}{
\partial x_{0}}\bigg|  \leq e^{N\kappa_{1}(t-t_{0})},
$$
where   $N$ depends only on $d$.
  Therefore, for any nonnegative
Lebesgue measurable function $w(x)$ it holds that
$$
e^{-N\kappa_{1}(t-t_{0})}
\int_{\bR^{d}}w(y)\,dy\leq
\int_{\bR^{d}}w(x_{t_{0},x_{0},t})\,dx_{0}\leq 
e^{N\kappa_{1}(t-t_{0})}\int_{\bR^{d}}w(y)\,dy .
$$
In particular, since
$$
\int_{\bR^{d}}|\chi^{(i)}_{t_{0},x_{0},t}( x)|^{2}\,dx_{0}=
\int_{\bR^{d}}|\chi^{(i)}(x-x_{t_{0},x_{0},t})|^{2}\,dx_{0} ,
$$
we have
$$
e^{-N\kappa_{1}(t-t_{0})}=N^{*}_{i} e^{-N\kappa_{1}(t-t_{0})}
 \int_{\bR^{d}}|\chi^{(i)}(x-y)|^{2}\,dy 
$$
$$
 \leq N^{*}_{i} 
\int_{\bR^{d}}|\chi^{(i)}_{t_{0},x_{0},t}( x)|^{2}\,dx_{0}
\leq N^{*}_{i} e^{N\kappa_{1}(t-t_{0})}
\int_{\bR^{d}}|\chi^{(i)}(x-y)|^{2}\,dy=e^{N\kappa_{1}(t-t_{0})}  ,
$$
where $N^{*}_{i}=|B_{1}|^{-1}
\rho_{0}^{-d}i^{d}$ and $|B_{1}|$ is the volume of  
$B_{1}$. It follows that
$$
\int_{\bR^{d}}|\chi^{(1)}_{t_{0},x_{0},t}( x)|^{2}\,dx_{0}
\leq (N^{*}_{1})^{-1}e^{N\kappa_{1}(t-t_{0})},
$$
$$
(N^{*}_{2})^{-1}e^{-N\kappa_{1}(t-t_{0})}\leq
\int_{\bR^{d}}|\chi^{(2)}_{t_{0},x_{0},t}( x)|^{2}\,dx_{0}.
$$
Furthermore, since $u_{t}=0$ if $\tau\geq t\geq t_{0}+\kappa^{-1}_{1}$
and $\chi^{(i)}_{t_{0},x_{0},t}=0$ if $t< t_{0}$,
in evaluating the norms in \eqref{3.14.2} we need not
integrate with respect to $t$ such that $\kappa_{1}(t-t_{0})\geq
1$, so that for all $t$ really involved we have
$$
\int_{\bR^{d}}|\chi^{(1)}_{t_{0},x_{0},t}( x)|^{2}\,dx_{0}
\leq (N^{*}_{1})^{-1}e^{N },\quad
(N^{*}_{2})^{-1}e^{-N }\leq
\int_{\bR^{d}}|\chi^{(2)}_{t_{0},x_{0},t}( x)|^{2}\,dx_{0}.
$$
After this observation it only remains to integrate 
\eqref{3.14.2} through with respect to $x_{0}$ and use the 
fact that $N^{*}_{1}=2^{-d}N^{*}_{2}$.
The lemma is proved.
 
{\bf Proof of Theorem \ref{theorem 3.11.1}}. First we show how to choose
$\gamma =\gamma (d,\delta,K)>0$. Call $N_{0}$
the constant factor of $\gamma ^{2/q} \| 
I_{t_{0}} Du\|_{\bL_{2}(\tau)}^{2}$
in \eqref{3.14.5}. We know that $N_{0}=N_{0}(d,\delta,K)$ and we choose
$\gamma \in(0,1]$ so that $N_{0}\gamma ^{2/q}\leq1/2$.
Then under the conditions of Lemma \ref{lemma 3.14.3}  
for $\lambda\geq1$
we have 
$$
\|I_{t_{0}} u\sqrt{c+\lambda }\|^{2}_{\bL _{2}(\tau)}+
\|I_{t_{0}} Du\|^{2}_{\bL _{2}(\tau)}
\leq N\big(\sum_{i=1}^{d}\| 
I_{t_{0}}f^{i}\|^{2}_{\bL _{2}(\tau)}
+\|I_{t_{0}} g\|^{2}_{\bL _{2}(\tau)}\big)
$$
$$
+N\lambda^{-1 }\|I_{t_{0}} f^{0}\|^{2}_{\bL _{2}(\tau)}
+NE\|u_{t_{0}}I_{t_{0}\leq\tau}\|^{2}_{\cL_{2}}
+
 N^{*} \lambda^{-1}\|  I_{t_{0}} Du\|_{\bL_{2}(\tau)}^{2}
$$
\begin{equation}
                                       \label{3.14.6}
+ N^{*}\| I_{t_{0}} u\|_{\bL_{2}(\tau)}^{2} 
+N^{*}\lambda^{-1}\sum_{i=1}^{d}\| 
I_{t_{0}}f^{i}\|^{2}_{\bL _{2}(\tau)}.
\end{equation}
After $\gamma $ has been fixed we have $\kappa_{1}=\kappa_{1}
(d,\delta,K,\rho_{0})$ and we take a $\zeta\in C^{\infty}_{0}(\bR)$
with support in $(0,\kappa_{1}^{-1})$
such that
\begin{equation}
                                       \label{3.15.1} 
\int_{-\infty}^{\infty}\zeta^{2}(t)\,dt=1.
\end{equation}
 For $s\in\bR$ define
$\zeta^{s}_{t}=\zeta(t-s)$, 
$u^{s}_{t}( x)=u_{t}(x)\zeta^{s}_{t}$. 
Obviously $u^{s}_{t}=0$ if $s_{+}+\kappa_{1}^{-1}\leq
t\leq\tau$.
Therefore,
we can apply 
\eqref{3.14.6}  to $u^{s}_{t}$ with $t_{0}=s_{+}$
 observing that
$$
du^{s}_{t}=(\sigma^{ik}_{t}D_{i}u^{s}_{t}+\nu^{}_{t}u^{s}_{t}
+\zeta^{s}_{t}g^{k})\,dw^{k}_{t}
$$
$$
+\big(D_{i}(a^{ij}_{t}D_{j}u^{s}_{t}
+\gb^{i}_{t}u^{s}_{t})+b^{i}_{t}u^{s}_{t}-(c_{t}+\lambda)
u^{s}_{t}+D_{i}(\zeta^{s}_{t}f^{i}_{t})+(
\zeta^{s}_{t}f^{0}_{t}+(\zeta^{s}_{t})'u_{t}\big)\,dt.
$$
Then from \eqref{3.14.6} for $\lambda\geq1$  we obtain
$$
  \|I_{s_{+}}\zeta^{s}  u\sqrt{c+\lambda }\|^{2}_{\bL _{2}(\tau)}+
\|I_{s_{+}}\zeta^{s}  Du\|^{2}_{\bL _{2}(\tau)}
$$
$$
\leq N\big(\sum_{i=1}^{d}\| I_{s_{+}}
\zeta^{s} f^{i}\|^{2}_{\bL _{2}(\tau)}
+\|I_{s_{+}}\zeta^{s}  g\|^{2}_{\bL _{2}(\tau)}\big)
$$
$$
+N\lambda^{-1 }\big(\|I_{s_{+}}\zeta^{s}f^{0}\|^{2}_{\bL _{2}(\tau)}
+\|I_{s_{+}}(\zeta^{s} )'u \|^{2}_{\bL _{2}(\tau)}\big)
+NE\|u_{s_{+}}\zeta^{s}_{s_{+}}I_{s_{+}\leq\tau}\|^{2}_{\cL_{2}}
$$
\begin{equation}
                                       \label{3.14.7} 
+ N^{*} \lambda^{-1}\|I_{s_{+}}\zeta^{s}Du\|_{\bL_{2}(\tau)}^{2}
+ N^{*}\| I_{s_{+}}\zeta^{s}  u\|_{\bL_{2}(\tau)}^{2}
+N^{*}\lambda^{-1}\sum_{i=1}^{d}\| 
I_{s_{+}}\zeta^{s} f^{i}\|^{2}_{\bL _{2}(\tau)}.
\end{equation}
Here $I_{s_{+}}$ can be dropped since 
$I_{s_{+}}I_{[0,\tau)}=I_{s}I_{[0,\tau)}$
and $I_{s}\zeta^{s}=\zeta^{s}$. After dropping
$I_{s_{+}}$ we 
 integrate through \eqref{3.14.7} with respect to $s\in\bR$, use
\eqref{3.15.1}, and observe that, since $\kappa_{1}$
depends only on $d,\delta,K,\rho_{0}$, we have
$$
\int_{-\infty}^{\infty}|\zeta'(s)|^{2}\,ds=N^{*}.
$$
We also use the fact that $\zeta^{s}_{s_{+}}\ne0$ only if
$s_{+}=0$ and
$-\kappa_{1}^{-1}\leq s\leq 0$ whereas
$$
\int_{-\kappa_{1}^{-1}}^{0}(\zeta^{s}_{0})^{2}\,ds=1.
$$

Then we conclude
$$
\lambda \| u\|^{2}_{\bL _{2}(\tau)}+\| u\sqrt{c}
\|^{2}_{\bL _{2}(\tau)}+
\| Du\|^{2}_{\bL _{2}(\tau)}
$$
$$
\leq N_{1}\big(\sum_{i=1}^{d}\| 
  f^{i}\|^{2}_{\bL _{2}(\tau)}
+\|  g\|^{2}_{\bL _{2}(\tau)}+E\|u_{0}\|^{2}_{\cL_{2}}\big)
$$
$$
+N_{1}\lambda^{-1 }\big(\|  f^{0}\|^{2}_{\bL _{2}(\tau)}
+\| u \|^{2}_{\bL _{2}(\tau)}\big)
+
 N_{1}^{*} \lambda^{-1}\|    Du\|_{\bL_{2}(\tau)}^{2}
$$
$$
+ N_{1}^{*}\|   u\|_{\bL_{2}(\tau)}^{2}
+N_{1}^{*}\lambda^{-1}\sum_{i=1}^{d}\| 
 f^{i}\|^{2}_{\bL _{2}(\tau)}.
$$
Without losing generality we assume that $N_{1}\geq1$
and  we show how to choose $\lambda_{0}=\lambda_{0}(
d,\delta,K,\rho_{0})$. We take it so that
  $\lambda_{0}\geq 4N^{*}_{1} $,
$\lambda_{0}^{2}\geq 4N_{1}$. Then we obviously come
to \eqref{3.11.2} with $N=4N_{1}$. The theorem is proved.

\mysection{Proof of Theorem \protect\ref{theorem 3.16.1}}
                                          \label{section 6.9.5}

We may assume in this section that $\cF_{t}=\cF_{t+}$
for all $t\in\bR_{+}$. This does not restrict generality
because replacing $\cF_{t}$ with $\cF_{t+}$ makes our 
assumptions weaker and does not affect our assertions
because the solutions are continuous in time.
Furthermore, having in mind setting all data equal to zero for
$t>\tau$, we see that without loss of generality we may assume
that $\tau=\infty$. Set   
$$
\bL_{2}=\bL_{2}( \infty),\quad
\bW^{1}_{2}=\bW^{1}_{2}( \infty),\quad
\cW^{1}_{2}=\cW^{1}_{2}( \infty).
$$
We need a few auxiliary results.
\begin{lemma}
                                         \label{lemma 6.28.1}
For any   $T,R\in\bR_{+}$, and $\omega\in\Omega$ we have  
\begin{equation}
                                                     \label{6.28.5}
\sup_{t\leq T}\int_{B_{R}}(|\gb_{t}(x)|^{q}
+|b_{t}(x)|^{q}+ c^{q}_{t}(x)   ) \,dx<\infty.
\end{equation}
\end{lemma}
Proof. Obviously it suffices to prove \eqref{6.28.5}
with $B_{\rho_{0}}(x_{0})$ in place of $B_{R}$ for any $x_{0}$.
In that case, for instance,  
$$
\int_{B_{\rho_{0}} (x_{0})} 
|\gb_{t}(x) |^{q}\,dx\leq 2^{q}
\int_{B_{\rho_{0}}(x_{0})} |\gb_{t}(x)-\bar{\gb}_{t}(x_{0})|^{q}\,dx
+N|\bar{\gb}_{t}(x_{0})|^{q}
$$
and we conclude estimating the left-hand side as in \eqref{6.11.2}
also relying on Assumption
\ref{assumption 3.16.1}.
Similarly, $b_{t}$ and $c_{t}$ are treated. The lemma is proved.

\begin{lemma}
                                         \label{lemma 3.16.1}
For any   $R\in\bR_{+}$
there exists a sequence of stopping times
$\tau_{n}\uparrow\infty$  such that
for any  $n=1,2,...$ and $\omega$ for almost any
$t\leq\tau_{n}$ we have
\begin{equation}
                                          \label{4.19.1}
\int_{B_{R}}(|\gb_{t}|^{q}+|b_{t}|^{q}+|c_{t}|^{q})
\,dx\leq n.
\end{equation}

\end{lemma}

Proof. For each $t,R>0$, and $\omega$ define
$$
\beta_{t,R}=  \int_{B_{R}}(|\gb_{t}|^{q}
+|b_{t}|^{q}+|c_{t}|^{q})\,dx,
$$
$$
\psi_{t,R}=\nlimsup_{\substack{0\leq s_{1}<s_{2}\leq t,\\
s_{2}-s_{1}\to0}}\frac{1}{s_{2}-s_{1}}\int_{s_{1}}^{s_{2}}
\beta_{s,R}\,ds.
$$
As is easy to see, $\psi_{t,R}$ is an increasing,
 left-continuous, and
$\cF_{t}$-adapted process. It follows that
$$
\tau_{n}:=\inf\{t\geq0:\psi_{t,R}> n\}
$$
are stopping times with respect to $\cF_{t+}$ ($=\cF_{t}$)
 and $\psi_{t,R}\leq n$ for $t<
\tau_{n}$. Furthermore, by Lemma \ref{lemma 6.28.1}
 we have
$\tau_{n}\uparrow\infty$ as $n\to\infty$. By Lebesgue 
differentiation theorem we conclude that
(for any $\omega$) for almost all $t\leq\tau_{n}$ we have
\eqref{4.19.1}. This proves the lemma.

By combining this lemma with Lemma \ref{lemma 6.27.1}
we obtain the following.

 \begin{corollary}
                                   \label{corollary 6.12.1}
If $\psi\in C^{\infty}_{0}$ has support in $B_{R}$,
then for $\tau_{n}$ from Lemma \ref{lemma 3.16.1}
for each    $n=1,2,...$, for almost
all $t\leq\tau_{n}$, for any $u\in W^{1}_{2}$ 
and $v\in W^{1}_{2} $
we have
$$
|(\gb^{i}_{t}D_{i}(v\psi),u)|\leq N\|v\|_{W^{1}_{2}}
\|u\|_{W^{1}_{2}},\quad
|(b^{i}_{t}D_{i}u,v\psi )|\leq N\|v\|_{W^{1}_{2}}
\|u\|_{W^{1}_{2}},
$$
\begin{equation}
                                             \label{6.12.4}
|(c_{t}v\psi,u)|\leq N\|v\|_{\cL_{2}}
\|u\|_{W^{1}_{2}},
\end{equation}
where the constant $N=N(n,d)$.
\end{corollary}

Since bounded linear operators
are continuous we obtain the following.
 \begin{corollary}
                                    \label{corollary 3.23.1}
If $\phi\in C^{\infty}_{0}$ has support in $B_{R}$,
then for $\tau_{n}$ from Lemma \ref{lemma 3.16.1}
and each $n$
 the operators
$$
u_{\cdot}\to (b^{i}_{\cdot}D_{i}u_{\cdot},\phi),\quad
u_{\cdot}\to (\gb^{i}_{\cdot}u_{\cdot},D_{i}\phi),
\quad
u_{\cdot}\to (c_{\cdot}u_{\cdot}, \phi),
$$
$$
u_{\cdot}\to \int_{0}^{\cdot}
(b^{i}_{t}D_{i}u_{t},\phi)\,dt,\quad
u_{\cdot}\to \int_{0}^{\cdot}
(\gb^{i}_{t}u_{t},D_{i}\phi)\,dt,\quad
u_{\cdot}\to \int_{0}^{\cdot}(c_{\cdot}u_{\cdot}, \phi)\,dt
$$
are continuous as operators from $\bW^{1}_{2}$ to
$\cL_{2}(\opar0,n\wedge\tau_{n}\cbrk)=
\cL_{2}(\opar0,n\wedge\tau_{n}\cbrk,\bR)$.
\end{corollary}

 In the proof of Theorem \ref{theorem 3.16.1} we are going to
use sequences which converge weakly in $\bW^{1}_{2}$.
Therefore, the following result is relevant.

\begin{lemma}
                                            \label{lemma 3.16.5}
Assume that for some
$f^{j}\in\bL_{2}$, $j=0,...,d$,
 $g=(g^{k})\in\bL_{2}$, $u\in\bW^{1}_{2}$,
and any $\phi\in C^{\infty}_{0}$
equation \eqref{3.16.7} with $u_{0}\in\cL_{2}(\Omega,\cF_{0},
\cL_{2})$ holds {\em
for almost all\/} $(\omega,t)$.
Then there exists a function $\tilde{u}\in\cW^{1}_{2}$
solving equation \eqref{2.6.4} (for all $t$)
with  initial data $u_{0}$
in the sense of Definition \ref{definition 3.20.01}.

\end{lemma}

Proof. We split the proof into two steps.

{\em Step 1. Modifying $u_{t}\psi$}.
We recall some facts from the theory of It\^o
stochastic integrals
in a separable 
Hilbert space, say $H$ and some other results, which
can be found, for instance, in \cite{KR} and \cite{Anal}. 
Integrating $H$-valued processes with respect to
a one-dimensional Wiener process presents no difficulties
and leads to strongly continuous $H$-valued locally
square-integrable martingales with natural isometry.
If $g=(g^{k})\in\bL_{2}$, then  
by Doob's inequality
$$
E\sup_{t}\big\|\sum_{k=n}^{m}
\int_{0}^{t} 
g^{k}_{s}\,dw^{k}_{s}\big\|^{2}_{\cL_{2} }
\leq 4E \int_{0}^{\infty} \sum_{k=n}^{m}\|
g^{k}_{s} \|^{2}_{\cL_{2} }\,ds\to0
$$
as $m\geq n\to\infty$. Therefore, 
$$
m_{t}=\sum_{k=1}^{\infty}\int_{0}^{t} 
g^{k}_{s}\,dw^{k}_{s} 
$$
is well defined as a continuous $\cL_{2}$-valued
square-integrable
martingale. Furthermore, for any $\phi\in C^{\infty}_{0}$
with probability one we have
$$
(m_{t},\phi)=\sum_{k=1}^{\infty}\int_{0}^{t}
(g^{k}_{s},\phi)\,dw^{k}_{s} 
$$
for all $t$ and the series on the right converges uniformly
in probability on $\bR_{+}$. If $g\in\bL_{2}
(\tau_{n})$, $n=1,2,...$, and stopping times $
\tau_{n}\uparrow\infty$, then  
$$
m_{t}=\sum_{k=1}^{\infty}\int_{0}^{t} 
g^{k}_{s}\,dw^{k}_{s} 
$$
is well defined as a locally square-integrable 
$\cL_{2}$-valued continuous martingale.
Again for any $\phi\in C^{\infty}_{0}$
with probability one we have
\begin{equation}
                                                    \label{3.23.1}
(m_{t},\phi)=\sum_{k=1}^{\infty}\int_{0}^{t}
(g^{k}_{s},\phi)\,dw^{k}_{s} 
\end{equation}
for all $t$ and the series on the right converges uniformly
in probability on every finite interval of time.

We fix a $\psi
\in C^{\infty}_{0}$ and apply the above to
$$
h^{\psi}_{t}:=\sum_{k=1}^{\infty}\int_{0}^{t} \psi(
\sigma^{ik}_{s}D_{i}u_{s}+\nu^{k}_{s}u_{s}+
g^{k}_{s})\,dw^{k}_{s}.
$$
Observe that, by assumption, for any $v\in C^{\infty}_{0}$
for almost all $(\omega,t) $
\begin{equation}
                                                      \label{3.17.3}
(u_{t}\psi,v)=(u_{0}\psi,v)+\int_{0}^{t}\langle F_{s},v\rangle\,ds
+(h^{\psi}_{t},v),
\end{equation}
where
$$
\langle F_{t},v\rangle=
 (b^{i}_{t}D_{i}u_{t}
-(c_{t}+\lambda)u_{t}+f^{0}_{t},v\psi)
-(a^{ij}_{t}D_{j}u_{t}+\gb^{i}_{t}u_{t}+
f^{i}_{t},D_{i}(v\psi)).
$$
We also define $V=W^{1}_{2}$,  and notice that if
$\|v\|_{V}\leq1$, then by Corollary \ref{corollary 6.12.1}
for any $T\in\bR_{+}$ for
almost any $(\omega,t)\in\Omega\times[0,T]$  we have
$$
|\langle F_{t},v\rangle|\leq N \big(
\sum_{j=0}^{d}\|f^{j}_{t}\|_{\cL_{2}}
+\|u_{t}\|_{W^{1}_{2}}\big),
$$
where $N$ is independent of $v,t$ 
(but may depend on $\omega$ and $T$).
It follows  that, for $V^{*}$
defined as the dual of $V$,
the $V^{*}$-norm of $F_{t}$ is in $\cL_{2}([0,T])$ (a.s.)
for every $T\in\bR_{+}$. It also follows that
\eqref{3.17.3} holds for almost all $(\omega,t)$
for each $v\in V$ rather than only for $v\in C^{\infty}_{0}$.

By Theorem 3.1 of \cite{KR} there exists a set $\Omega_{\psi}$
of full probability and an $\cL_{2}$-valued
function
$\tilde{u}^{\psi}_{t}$ on $\Omega\times\bR_{+}$ such that
$\tilde{u}^{\psi}_{t} $ is $\cF_{t}$-measurable,
$\tilde{u}^{\psi}_{t}$ is $\cL_{2}$-continuous in $t$
 for every $\omega$ and 
$\tilde{u}^{\psi}_{t}=u_{t}\psi$ for almost
all $(\omega,t)$. 
Furthermore,
for $\omega\in\Omega_{\psi}$, $t\geq0$,
and $\phi\in C^{\infty}_{0}$ we have
$$
(\tilde{u}^{\psi}_{t},\phi)=(h^{\psi}_{t},\phi)
+\int_{0}^{t} (b^{i}_{s}D_{i}u_{s}
-(c_{s}+\lambda)u_{s}+f^{0}_{s},\phi\psi)\,ds
$$
\begin{equation}
                                                  \label{3.19.3}
-\int_{0}^{t} \big(a^{ij}_{s}D_{j}u_{s}+\gb^{i}_{s}u_{s}+
f^{i}_{s},D_{i}(\phi\psi)\big)\,ds.
\end{equation}

{\em Step 2. Constructing $\tilde{u}_{t}$\/}.
Let $\psi\in C^{\infty}_{0}$ be such that
$\psi=1$ on $B_{1}$ and set $\psi_{n}(x)=\psi(x/n)$,
$n=1,2,...$. Define $\tilde{u}^{n}_{t}=
\tilde{u}^{\psi_{n}}_{t}$ and notice that by the above
 for $m\geq n$
and almost all $(\omega,t)$
$$
\tilde{u}^{m}_{t}I_{B_{n}}=u_{t}\psi_{m}I_{B_{n}}
=u_{t}I_{B_{n}}=\tilde{u}^{n}_{t}I_{B_{n}}
 $$
 as $\cL_{2}$-elements.
Since the extreme terms
 are $\cL_{2}$-continuous functions of $t$,
there exist sets $\Omega_{nm}$, $m\geq n$, of full probability
such that for $\omega\in\Omega_{nm}$ we have
$\tilde{u}^{m}_{t}I_{B_{n}}=\tilde{u}^{n}_{t}
I_{B_{n}}$ as $\cL_{2}$-elements 
for   all $t$.

Then for $t\geq0$
and  $\omega\in\Omega':=\bigcap_{m\geq n}\Omega_{nm}$ the formula
$$
\tilde{u}_{t}=I_{\Omega'}
\sum_{n=0}^{\infty}\tilde{u}_{t}^{n+1}
I_{B_{n+1}\setminus B_{n}}
$$
defines a distribution such that $\tilde{u}_{t}I_{B_{n}}=
\tilde{u}^{n}_{t}I_{B_{n}}$ as $\cL_{2}$-elements
for any $\omega\in\Omega'$, $t\geq0$, and $n$. 
It follows that $\tilde{u}_{t}=u_{t}$
 as distributions for almost any $(\omega,t)$, hence,
$\tilde{u}\in\bW^{1}_{2}$ and there exists an event
$\Omega''\subset\Omega'$ of full probability such that
for any $\omega\in\Omega''$ and almost any $t\geq0$
we have $\tilde{u}_{t}=u_{t}$. Now
\eqref{3.19.3} implies that if $\phi\in C^{\infty}_{0}$
is such that $\phi(x)=0$ for $|x|\geq n$, then
for $\omega\in\Omega''\cap\Omega_{\psi_{n}}$
and all $t\geq0$ we have
$$
(\tilde{u}_{t},\phi)=(\tilde{u}^{n}_{t},\phi)
 =(h^{\psi_{n}}_{t},\phi)
+\int_{0}^{t} (b^{i}_{s}D_{i}\tilde{u}_{s}
-(c_{s}+\lambda)\tilde{u}_{s}+f^{0}_{s},\phi)\,ds
$$
\begin{equation}
                                                  \label{3.19.03}
-\int_{0}^{t} \big(a^{ij}_{s}D_{j}\tilde{u}_{s}+
\gb^{i}_{s}\tilde{u}_{s}+
f^{i}_{s},D_{i}\phi \big)\,ds.
\end{equation}

By recalling what was said about \eqref{3.23.1}
and using Corollary \ref{corollary 6.12.1}, we see that
indeed the requirements   of Definition
\ref{definition 3.20.01} are satisfied 
 with $\tilde{u}$ and $\infty$ in place of
$u$ and $\tau$, respectively. The lemma is proved.

\begin{lemma}
                                            \label{lemma 3.23.2}
Let $\phi\in C^{\infty}_{0}$ be supported in $B_{R}$ and take
  $\tau_{n}$ from Lemma \ref{lemma 3.16.1}.
Let $u^{n}$, $u\in \bW^{1}_{2}$, $n=1,2,...$, be such that
$u^{n}\to u$ weakly in $\bW^{1}_{2}$. 
For $n=1,2,...$   define $\chi_{n}(t)=(-n)\vee t\wedge n$,
$\gb^{i}_{nt}=\chi_{n}(\gb^{i}_{t})$, 
$b^{i}_{nt}=\chi_{n}(b^{i}_{t})$  and 
set $c_{ns}=n\wedge c_{s}$.
Then for any   $m=1,2,...$
$$
\int_{0}^{t}[(b^{i}_{ns}D_{i}u^{n}_{s},\phi)
-(\gb^{i}_{ns}u^{n}_{s},D_{i}\phi)-(c_{ns}u^{n}_{s},\phi)]\,ds
$$
\begin{equation}
                                                 \label{4.19.6}
\to
\int_{0}^{t}[(b^{i}_{s}D_{i}u_{s},\phi)
-(\gb^{i}_{s}u_{s},D_{i}\phi)-(c_{ s}u _{s},\phi)]\,ds
\end{equation}
weakly in the space $\cL_{2}(\opar0,m\wedge\tau_{m}\cbrk)$
as $n\to\infty$ .
\end{lemma}

Proof. By Corollary \ref{corollary 3.23.1} and by the
fact that (strongly) continuous operators are weakly
continuous we obtain that
$$
\int_{0}^{t}[(b^{i}_{s}D_{i}u^{n}_{s},\phi)
-(\gb^{i}_{s}u^{n}_{s},D_{i}\phi)-(c_{ s}u^{n}_{s},\phi)]\,ds
$$
$$
\to
\int_{0}^{t}[(b^{i}_{s}D_{i}u_{s},\phi)
-(\gb^{i}_{s}u_{s},D_{i}\phi)-(c_{ s}u _{s},\phi)]\,ds
$$
as $n\to\infty$ 
weakly in the space $\cL_{2}(\opar0,m\wedge\tau_{m}\cbrk)$
for any $m$. Therefore, it suffices to show that
$$
\int_{0}^{t}[(D_{i}u^{n}_{s},(b^{i}_{s}-b^{i}_{ns})\phi)
-(u^{n}_{s},(\gb^{i}_{s}-\gb^{i}_{ns})D_{i}\phi
+(c_{s}-c_{ns})\phi)]\,ds
\to0
$$
weakly in $\cL_{2}(\opar0,m\wedge\tau_{m}\cbrk)$
for any $m$. In other words, it suffices to show
that for any $\xi\in \cL_{2}(\opar0,
m\wedge\tau_{m}\cbrk)$
$$
E\int_{0}^{m\wedge\tau_{m}}\xi_{t}
\big(\int_{0}^{t}[(D_{i}u^{n}_{s},(b^{i}_{s}-b^{i}_{ns})\phi)
$$
$$
-(u^{n}_{s},(\gb^{i}_{s}-\gb^{i}_{ns})D_{i}\phi 
+(c_{s}-c_{ns})\phi)]\,ds
\big)\,dt\to0.
$$
This relation is rewritten as
$$
E\int_{0}^{m\wedge\tau_{m}}
[(D_{i}u^{n}_{s},\eta_{s}(b^{i}_{s}-b^{i}_{ns})\phi)
$$
\begin{equation}
                                                 \label{4.21.1}
-(\eta_{s}u^{n}_{s},
(\gb^{i}_{s}-\gb^{i}_{ns})D_{i}\phi
+(c_{s}-c_{ns})\phi)]\,ds\to0,
\end{equation}
where the process
$$
\eta_{s}:=\int_{s}^{m\wedge\tau_{m}}\xi_{t}\,dt
$$
is of class $\cL_{2}(\opar0,m\wedge\tau_{m}\cbrk)$ since
$m\wedge\tau_{m}$ is bounded ($\leq m$).

However, by the choice of $\tau_{m}$ and the 
dominated convergence theorem,
$$
\eta_{s}(\gb^{i}_{s}-\gb^{i}_{ns})D_{i}\phi\to0,\quad
\eta_{s}(b^{i}_{s}-b^{i}_{ns})\phi\to0,\quad
\eta_{s}(c_{s}-c_{ns})\phi\to0
$$
as $n\to\infty$
strongly in $\bL_{2}(\opar0,m\wedge\tau_{m}\cbrk)$ 
(use   the fact that $q\geq2$) and 
by assumption $u^{ n}\to u$ and $Du^{n}\to Du$
weakly in $\bL_{2}(\opar0, \tau_{m}\cbrk)$. This implies
\eqref{4.21.1}  for any $m$  and the lemma is proved.

 {\bf Proof of Theorem \ref{theorem 3.16.1}}. Define
$\gb_{nt}$, $b_{nt}$, and $c_{nt}$  as in  
  Lemma \ref{lemma 3.23.2} and consider equation
\eqref{2.6.4} with $\gb_{nt}$, $b_{nt}$, and $c_{nt}$ 
in place of $\gb_{t}$, $b_{t}$, and $c_{t}$, respectively,
and with $\tau=n$. By a classical result there exists
a unique $u^{n}\in\cW^{1}_{2}(n)$ satisfying the modified
equation with initial condition $u_{0}$. 
Obviously, $\gb_{nt}$, $b_{nt}$, and $c_{nt}$  satisfy Assumption
\ref{assumption 3.11.1} with the same $\gamma $
as $\gb_{t}$, $b_{t}$, and $c_{t}$ do.
By Theorem \ref{theorem 3.11.1} for
$\lambda\geq\lambda_{0}(d,\delta,K,\rho_{0})$ we have
$$
\|u^{n}\|_{\bL_{2}(n)}+\|Du^{n}\|_{\bL_{2}(n)}\leq N,
$$
where $N$ is independent of $n$. Hence the sequence of functions
$u^{n}_{t}I_{t\leq n}$ is bounded in the Hilbert space 
$\bW^{1}_{2}$ and consequently has a weak limit 
point $u\in \bW^{1}_{2}$. For simplicity of presentation
we assume that the whole sequence $u^{n}_{t}I_{t\leq n}$
converges weakly to $u$. 
Take a $\phi\in C^{\infty}_{0}$. Then
by Lemma \ref{lemma 3.23.2} for appropriate $\tau_{m}$ we have
that \eqref{4.19.6} holds weakly in
$\cL_{2}(\opar0,m\wedge\tau_{m}\cbrk)$
for any $m$. Since
$$
u=u_{t}\to\sum_{k=1}^{\infty}
\int_{0}^{t}(\Lambda^{k}_{s}u_{s},\phi)\,dw^{k}_{s}
$$
is a continuous operator from $\bW^{1}_{2}$ to
$\bL_{2}(\opar0,m\cbrk)$, it is weakly continuous, so that
$$
\sum_{k=1}^{\infty}
\int_{0}^{t}(\Lambda^{k}_{s}u^{n}_{s},\phi)\,dw^{k}_{s}
\to \sum_{k=1}^{\infty}
\int_{0}^{t}(\Lambda^{k}_{s}u_{s},\phi)\,dw^{k}_{s}
$$
weakly in $\cL_{2}(\opar0,m\cbrk)$ for any $m$. Obviously,
the same is true for $(u^{n}_{ t},\phi)\to(u_{t},\phi)$
and the remaining terms entering 
 the equation for  $ u^{n}_{ s}$.
Hence by passing to the weak limit in the equation
for $u^{n}_{ t}$ we see that $u$ satisfies the assumptions
of Lemma \ref{lemma 3.16.5} applying which
finishes the proof of the theorem.


\begin{thebibliography}{mm}

 \bibitem{AM} S. Assing
and R. Manthey, {\em Invariant measures for
stochastic heat equations with unbounded coefficients\/}, 
Stochastic Process. Appl.,  Vol. 103  (2003),  No. 2, 237-256.
 
 \bibitem{CV1} P. Cannarsa  and V. Vespri, {\em Generation
of analytic semigroups by elliptic operators with unbounded
coefficients\/}, SIAM J. Math. Anal., Vol. 18 (1987),
No. 3, 857-872.

\bibitem{CV} P. Cannarsa  and V. Vespri, {\em
 Existence and uniqueness results
for a nonlinear stochastic partial differential equation\/}, in 
Stochastic Partial Differential Equations and Applications
Proceedings, G. Da Prato and L. Tubaro (eds.), Lecture Notes in
Math., Vol. 1236, pp. 1-24, Springer Verlag, 1987.

\bibitem{ChG} A. Chojnowska-Michalik and B. Goldys,  {\em
Generalized symmetric Ornstein-Uhlenbeck semigroups in $L^p$:
Littlewood-Paley-Stein inequalities and domains of
generators\/}, J. Funct. Anal., Vol. 182 (2001), 243-279.


\bibitem{CF} G. Cupini and S. Fornaro, {\em Maximal regularity
in $L^ p(\bR^N)$ for a class of elliptic operators with unbounded
coefficients\/},  Differential Integral Equations, Vol.  17 
(2004), 
 No. 3-4, 259-296.



 \bibitem{GL} M. Geissert and A. Lunardi,  {\em Invariant
measures and maximal $L\sp 2$ regularity for nonautonomous
Ornstein-Uhlenbeck equations\/},  J. Lond. Math. Soc. (2), Vol.
  77  (2008), No. 3, 719-740.

\bibitem{FL} B. Farkas and A. Lunardi, {\em Maximal regularity for Kolmogorov
operators in $L^2$ spaces with respect to invariant measures\/}, J.
Math. Pures Appl., Vol. 86 (2006), 310-321.

\bibitem{Gy93}   I. Gy\"ongy,
{\em Stochastic partial differential equations on
Manifolds, I\/}, Potential Analysis, Vol.  2 (1993), 101-113.

\bibitem{Gy97}   I. Gy\"ongy, {\em
Stochastic partial differential equations
manifolds II. Nonlinear filtering\/},
Potential Analysis, Vol. 6 (1997),  39-56.

\bibitem{GK}   I. Gy\"ongy and N.V. Krylov,  {\em
  On stochastic partial differential equations with unbounded
 coefficients\/},  Potential Analysis, Vol. 1 (1992), No. 3, 233-256.

\bibitem{Ki1} Kyeong-Hun Kim,  {\em
 On $L_p$-theory of stochastic partial differential equations of
divergence form in $C^1$ domains\/}, Probab. Theory Related Fields,
Vol. 130 (2004), No. 4, 473-492.

\bibitem{Anal} N.V. Krylov,
 {\em An analytic approach to SPDEs}, pp. 185-242 in
Stochastic Partial Differential Equations: Six Perspectives,
Mathematical Surveys and Monographs, Vol. 64,
AMS, Providence, RI, 1999.

 \bibitem{Kr09_1} N.V. Krylov,
 {\em On   linear
elliptic and parabolic equations
with growing drift in Sobolev spaces without weights\/},
Problemy Matemtaticheskogo Analiza, Vol. 40 (2009), 
77-90, in Russian;
English version in Journal of 
Mathematical Sciences, Vol. 159 (2009), No. 1, 75-90, Srpinger.

\bibitem{Kr09_2} N.V. Krylov,
 {\em
On divergence
form SPDEs with VMO coefficients\/},   
SIAM J. Math. Anal. Vol. 40 (2009), No. 6,  2262-2285.

\bibitem{Kr09_3} N.V. Krylov,
 {\em It\^o's formula
for the $L_{p}$-norm of
  stochastic $W^{1}_{p}$-valued processes\/}, to appear in 
Probab. Theory Related Fields,
 http://arxiv.org/abs/0806.1557 

\bibitem{Kr09_4} N.V. Krylov,
{\em On  the It\^o-Wentzell formula for
distribution-valued processes and related topics\/},
submitted to Probab. Theory Related Fields,
 http://arxiv.org/abs/0904.2752 

 

\bibitem{Kr_10} N.V. Krylov, {\em Filtering equations for partially
observable diffusion processes with Lipschitz
continuous coefficients\/}, to appear in
``The Oxford
Handbook of Nonlinear Filtering", Oxford University Press,
http://arxiv.org/abs/0908.1935

\bibitem{KP} N.V. Krylov and E. Priola, {\em
Elliptic and parabolic second-order PDEs with growing coefficients\/},
to appear in Comm. in PDEs, http://arXiv.org/abs/0806.3100
  
\bibitem{KR} N.V. Krylov and
  B.L. Rozovskii,  {\em  Stochastic evolution
equations\/}, pp.  71-146 in ``Itogy nauki i
tekhniki'',  Vol. 14, VINITI, Moscow, 1979, in Russian;
English translation: J. Soviet Math.,
 Vol. 16 (1981), No. 4, 1233-1277.



\bibitem{Lu}
A. Lunardi, {\em Schauder estimates for a class of degenerate
elliptic and parabolic operators with unbounded coefficients in
$\bR^n$\/}, Ann. Sc. Norm. Super Pisa, Ser. IV., Vol. 24 (1997), 133Ð164.


\bibitem{LV} A. Lunardi and V. Vespri,
{\em Generation of strongly continuous semigroups by elliptic
operators with unbounded coefficients in $L^p(\bR^n)$\/},
Rend. Istit. Mat. Univ. Trieste 28 (1996), suppl., 
251-279 (1997).

\bibitem{MP1}
G. Metafune, J. Pr\"uss, A. Rhandi, and R. Schnaubelt, {\em The domain
of the OrnsteinÐUhlenbeck operator on an $L^p$-space with invariant
measure\/}, Ann. Sc. Norm. Super. Pisa, Cl. Sci., (5) 1 (2002),
471-485.

\bibitem{MP} G.
Metafune, J. Pr\"uss,  A. Rhandi, and R. Schnaubelt,  {\em
$L\sp p$-regularity for elliptic operators with unbounded 
coefficients\/},
Adv. Differential Equations, Vol. 10 (2005), No. 10, 1131-1164. 

\bibitem{PR} J.
Pr\"uss, A. Rhandi,  and R. Schnaubelt, {\em
The domain of elliptic operators on $L^p(\bR^d)$ with
unbounded drift coefficients\/}, Houston J. Math., Vol. 32
(2006), No. 2, 563-576.

 
\end{thebibliography}
\end{document}